\documentclass [12pt,frenchy]{article}
\usepackage{fancyhdr,amsmath, graphicx, psfrag, color ,amsfonts, layout, enumerate, lscape, tikz}
\usetikzlibrary{shapes}

\font\fontauthors=cmcsc10 scaled \magstep1

\def\be{\begin{equation}}

\def\conv{\smash{\mathop{\longrightarrow}\limits _{n\to \infty}}}

\def\convloi{\smash{\mathop{~~\longrightarrow~}\limits _{n\to \infty}\limits ^{\cal D}}}

\def\doubleindice#1#2{\genfrac{}{}{0pt}{1}{#1}{#2}}
\def\ee{\end{equation}}
\def\egalenloi{\smash{\mathop{~= ~}\limits ^{\cal L}}}
\def\egalLoi{{~\mathop{= }\limits^{\rond L}}~}

\def\summ#1#2{\sum _{\begin{array}{cc}\scriptstyle #1\\\scriptstyle #2\end{array}}}

\def\Supp{\mathop{\rm Supp}\nolimits}
\def\stirling2 #1#2{\left\{\begin{matrix} #1\\#2\end{matrix}\right\}}
\def\summ#1#2{{\displaystyle\sum _{\doubleindice{#1}{#2}}}}

\def\Var{\mathop{\rm Var}\nolimits}

\newcommand{\pff}{\noindent {\sc Proof.}\ }

\def\QED{\hfill\vrule height 1.5ex width 1.4ex depth -.1ex \vskip 10pt}
\newcommand\1{\leavevmode\hbox{\rm \small1\kern-0.35em\normalsize1}}

\def\g#1{\mathbb #1}
\def\rond#1{\mathcal #1}

\def\E{\g E}

\newtheorem {Th}{Theorem}
\newtheorem {Cor}{Corollary}
\newtheorem {Prop}{Proposition}
\newtheorem {Lem}{Lemma}
\newtheorem {Rem}{Remark}

\begin{document}

\begin{center}
\LARGE{\bf {Smoothing equations for large P\'olya urns\footnote{

Urn model.
Martingale.
Multitype branching process.
Smoothing transforms.
Contraction method.
Characteristic function.
Moment-determined probability distributions.
}}}\\
\end{center}

\begin{center}
May 29th 2013\\
%by
\end{center}

\begin{center}
{\fontauthors
Brigitte Chauvin,
%\footnote{INRIA Rocquencourt, project Algorithms and Laboratoire de Math\'ematiques de Versailles, CNRS, UMR 8100 --  INRIA Domaine de Voluceau B.P.105, 78153 Le Chesnay CEDEX (France).},
%Quansheng Liu,
C\'ecile Mailler,
Nicolas Pouyanne,
}
\medskip
%\footnote{Laboratoire de Math\'ematiques de Versailles, CNRS, UMR 8100 -- Universit\'e de Versailles - St-Quentin,  45 avenue des Etats-Unis, 78035 Versailles CEDEX (France).
\\
Universit\'e de Versailles-St-Quentin,\\
Laboratoire de Math\'ematiques de Versailles,\\
CNRS, UMR 8100,\\
45, avenue des Etats-Unis, \\
78035 Versailles CEDEX, France.
\end{center}

\vskip 1truecm
\noindent{\bf Abstract.}
Consider a balanced non triangular two-color P\'olya-Eggenberger urn process, assumed to be large which means that the ratio $\sigma$ of the replacement matrix eigenvalues satisfies $1/2<\sigma <1$.
The composition vector of both discrete time and continuous time models admits a drift which is carried by the principal direction of the replacement matrix.
In the second principal direction, this random vector admits also an almost sure asymptotics and a real-valued limit random variable arises, named $W^{DT}$ in discrete time and $W^{CT}$ in continous time.
The paper deals with the distributions of both $W$.
Appearing as martingale limits, known to be nonnormal, these laws remain up to now rather mysterious.

Exploiting the underlying tree structure of the urn process, we show that $W^{DT}$ and $W^{CT}$ are the unique solutions of two distributional systems in some suitable spaces of integrable probability measures.
These systems are natural extensions of distributional equations that already appeared in famous algorithmical problems like Quicksort analysis.
Existence and unicity of the solutions of the systems are obtained by means of contracting smoothing transforms.
Via the equation systems, we find upperbounds for the moments of $W^{DT}$ and $W^{CT}$
and we show that the laws of $W^{DT}$ and $W^{CT}$ are moment-determined.
We also prove that $W^{DT}$ is supported by the whole real line and admits a continuous density
($W^{CT}$ was already known to have a density, infinitely differentiable on $\g R\setminus\{ 0\}$ and not bounded at the origin).

\tableofcontents

%%%%%%%%%%%%%%%%%%%%
\section{Introduction}
\label{intro}
%%%%%%%%%%%%%%%%%%%%%

P\'olya urns provide a rich model for many situations in algorithmics.
Consider an urn that contains red and
black balls. 
Start with a finite number of red and black balls as initial composition (possibly monochromatic).
At each discrete time $n$,
 draw a ball at random, notice its color,
put it back into the urn and add balls according to the following rule:
if the drawn ball  is red,  add $a$ red balls and $b$ black balls;
if the drawn ball  is black, add $c$ red balls and $d$ black balls.
The integers $a,b,c,d$ are assumed to be nonnegative\footnote
{One admits classically negative values for $a$ and $d$, together with arithmetical conditions on $c$ and $b$.
Nevertheless, the paper deals with so-called \emph{large} urns, for which this never happens.}.
Thus, the replacement rule is described by the so-called \emph{replacement matrix} $$
R=\left(
\begin{array}{cc}
a&b\\c&d
\end{array}
\right).
$$
``Drawing a ball at random'' means choosing \emph{uniformly} among the balls contained in the urn.
That is why this model is related to many situations in mathematics, algorithmics or theoretical physics where a uniform choice among objects determines the evolution of a process. See Johnson and Kotz's book \cite{JK}, Mahmoud's book \cite{Mah08} or Flajolet et al.  \cite{FlaDumPuy} for many examples.

In the present paper, the urn is assumed to be {\it balanced}, which means that the total number of balls added at each
step is a constant
$$
S=a+b=c+d.
$$
The composition vector of the urn at time $n$ is denoted by
$$
U^{DT}(n)=\left(
\begin{array}{c}
\hbox{number of red balls at time}~n\\
\hbox{number of black balls at time}~n
\end{array}
\right).
$$
Two main points of view are classically used on this random vector.
The \emph{forward} point of view consists in considering the composition vector sequence $\left(U^{DT}(n)\right) _{n\in\g N}$ as a $\g N^2$-valued Markov chain.
%This follows from the fact that the random composition at a given time depends only on the probability distribution of the preceding composition.
The information on the successive configurations is thus concentrated in a global object:
the random process, giving access to probabilistic tools like martingales, embedding in continuous time, branching processes.
A vast part of the literature on P\'olya urns relies on such probability tools, dealing most often with natural extensions of the model to a random replacement matrix or to an arbitrary finite number of colors.
The forward point of view is particularly efficient to get results on the asymptotics of the process.
See for instance Janson's seminal paper \cite{Jan} or \cite{Pou08} for an extensive state of the art on such methods.

Alternatively, a natural feature consists in using the recursive properties of the random structure through a {\emph{divide and conquer}} principle.
This is the \emph{backward} point of view.
Applied to generating functions, it is the base tool for analytic combinatorics methods, developed in Flajolet \emph{et al.} papers \cite{FlaGabPek, FlaDumPuy}.
Expressed in terms of the random process, the backward approach leads to dislocation equations on limit distributions that can already be found in a wide generality in Janson \cite{Jan};
these equations are further developed in \cite{ChaPouSah} for two-color urns and in \cite{ChaLiuPouContinu,ChaLiuPouDiscret} for the urn related to $m$-ary search trees as well.

\vskip 5pt
In order to state our results and also the asymptotic theorems they are based on, we first give some notations that are made more complete in Section~\ref{sec:model}.
The eigenvalues of the replacement matrix $R$ are $S$ and the integer
$$
m:= a-c = d-b
$$
and we denote by
$$
\sigma := \frac mS\leq 1
$$
the ratio between these eigenvalues.
The particular case $\sigma =1$ is the original P\'olya urn (see P\'olya~\cite{Polya});
this process has a specific well known asymptotics with a random drift.
In appendix, our Section~\ref{sec:appendix} is devoted to gather results on this almost sure limit and on the asymptotic Dirichlet distribution as well.
When $\sigma <1$, it is well known that the asymptotics of the process has two different behaviours, depending on the position of $\sigma$ with respect to the value $1/2$
(see Athreya and Karlin \cite{AK} for the original result, Janson~\cite{Jan} or~\cite{Pou08} for the results below).
Briefly said,

\vskip 5pt
\noindent
{\it (i)}
when $\sigma <\frac{1}{2} $, the urn is called {\bf\emph{small}} and, except when $R$ is triangular,
the composition vector is asymptotically Gaussian\footnote
{The case $\sigma=1/2$ is similar to this one, the normalisation being $\sqrt{n\log n}$ instead of
$\sqrt n$.}:
$$
\frac{U^{DT}\left( n\right) -nv_1}{\sqrt n} \convloi\rond G\left( 0,\Sigma^2\right)
$$
where $v_1$ is a suitable eigenvector of $^t\! R$ relative to $S$ and $\rond G$ a centered Gaussian vector with covariance matrix $\Sigma^2$ that has a simple
closed form;

\vskip 5pt
\noindent
{\it (ii)}
when $\frac{1}{2}  < \sigma  < 1  $, the urn is  called {\bf\emph{large}} and the composition vector
has a quite different strong asymptotic form:
\begin{equation}
\label{asymptotiqueDT}
U^{DT}\left( n\right) = nv_1 + n^{\sigma} W^{DT} v_2 + o\left( n^{\sigma}\right)
\end{equation}
where 
$v_1, v_2$ are suitable (non random) eigenvectors of $^t\! R$ relative to the respective eigenvalues $S$ and $m$,
$W^{DT}$ is a real-valued random variable arising as the limit of a martingale,
the little $o$ being almost sure and in any $L^p, p\geq 1$.
%The moments of $W^{DT}$ can be recursively calculated but they have no known global closed form \cite{Pou08}.

\vskip 5pt
Classically, like for any Markov chain, one can embed the discrete time process $\left( U^{DT}(n)\right) _{n\in\g Z_{\geq 0}}$ into continuous time.
In the case of P\'olya urns having a replacement matrix with nonnegative entries, this defines a two-type branching process
$$
\left( U^{CT}(t)\right) _{t\in\g R_{\geq 0}}.
$$
A similar phase transition occurs when $t$ tends to infinity:
for small urns, the process $U^{CT}$ has a (random) almost sure drift and satisfies a gaussian central limit theorem
(see Janson~\cite{Jan}).
When the urn is large, the asymptotic behaviour of the process, when $t$ tends to infinity, is given by 
$$
U^{CT}\left( t\right) =e^{St} \xi v_1\left( 1+o(1)\right) +e^{mt}W^{CT}v_2\left( 1+o(1)\right) ,
$$
where $\xi$ is {\rm Gamma}-distributed, $W^{CT}$ is a real-valued random variable arising as the limit of a martingale,
the little $o$ is almost sure and in any $L^p, p\geq 1$, the basis $\left( v_1,v_2\right)$ of deterministic vectors being the same one as in~\eqref{asymptotiqueDT}.
These asymptotic results are more detailed in Section~\ref{sec:model}.
Because of the canonical link between $U^{DT}$ and $U^{CT}$ \emph{via} stopping times, the two random variables $W^{DT}$ and $W^{CT}$ are related by the so-called \emph{martingale connexion} as explained in Section~\ref{sec:martingaleConnexion}.
Consequently any information about one distribution is of interest for the other one.
All along the paper, the symbol $DT$ is used to qualify \emph{discrete-time} objects while $CT$ will refer to \emph{continuous-time} ones.

\vskip 5pt
In this article, we are interested by large urns.
More precisely, the attention is focused on the non classical distributions in $W^{DT}$ and $W^{CT}$ when the replacement matrix $R$ is not triangular (\emph{i.e.} when $bc\neq 0$).
For example, $W^{CT}$ is not normally distributed, which can be seen on its exponential moment generating series that has a radius of convergence equal to zero, as shown in~\cite{ChaPouSah}
(see Section~\ref{sec:moments} for more details).
Because of the martingale connexion, this implies that $W^{DT}$ is not normal either.
Our main goal is to get descriptions of these laws (density, moments, tail, \dots).

What is already known about $W^{DT}$ or $W^{CT}$?
In \cite{ChaPouSah}, the Fourier transform of $W^{CT}$ is ``explicitely'' calculated, in terms of the inverse of an abelian integral on the Fermat curve of degree $m$.
%It was also proved that the Laplace series of the moments has a radius of convergence equal to zero.
The existence of a density with respect to the Lebesgue measure on $\g R$ and the fact that $W^{CT}$ is supported by the whole real line are deduced from this closed form.
Nevertheless, the order of magnitude of the moments and the question of the determination of the law by its moments remained open questions. The shape of the density was mysterious, too. The present paper answers to these questions in Section \ref{sec:moments} and \ref{sec:shape} respectively.

In the present text, we exploit the underlying tree structure of a P\'olya urn.
Governing both the backward and the forward points of view, it contains a richer structure than the plain composition vector process.
Section \ref{sec:decomposition} is devoted to highlighting this tree process and to 
derive decomposition properties on the laws of the composition vector at finite time.
These decompositions directly lead to distributional fixed point systems~\eqref{pointFixeDiscret} and~\eqref{pointFixeContinu} respectively satisfied by $W^{DT}$ and $W^{CT}$, as stated in Theorem~\ref{th:systemDT} and Theorem~\ref{th:systemCT}.

With a slightly different approach, Knape and Neininger \cite{KnaNei} start from the tree decomposition of the discrete P\'olya urn and establish the fixed point system \eqref{pointFixeDiscret} with the contraction method tools developed in Neininger-R\"uschendorf \cite{NeiRusaap}. This complementary point of view does not take advantage of the limit random variable $W^{DT}$ but applies for small and large urns together, allowing to find limit Gaussian distributions thus providing an alternative method to the embedding method used by Janson in \cite{Jan}.

Sometimes called fixed point equations for the smoothing transform or just \emph{smoothing equations} in the literature (Liu \cite{Liu01}, Durrett-Liggett \cite{DL83})
), distributional equations of type
\begin{equation}
\label{smoothingEq}
X\egalenloi \sum A_iX^{(i)}
\end{equation}
have given rise to considerable interest in, and literature on. For a survey, see Aldous-Bandyopadhyay \cite{AldBan}. 
In theoretical probability, they are of relevance in connexion with branching processes (like in Liu \cite{Liu98}, Biggins-Kyprianou  \cite{BigKyp05}, Alsmeyer et al \cite{AlsBigMai}) or with Mandelbrot cascades (Mandelbrot \cite{Man}, Barral \cite{Bar}). 
They occur in various areas of applied probability, 
and also on the occasion of famous problems arising in analysis of algorithms, like Quicksort (R\"osler \cite{Ros92}). They are naturally linked with the analysis of recursive algorithms and data structures 
(Neininger-R\"uschendorf \cite{NeiRus}, surveys in Rosler-R\"uschendorf \cite{RR} or Neininger-R\"uschendorf \cite{NeiRusSurvey})

Most often, in Equation~\eqref{smoothingEq}, the $A_i$ are given random variables and the $X^{(i)}$ are independent copies of $X$, independent of the $A_i$ as well.
Our System~\eqref{pointFixeContinu} with unknown real-valued random variables (or distributions) $X$ and $Y$ is the following:
$$
\left\{
\begin{array}{l}
\displaystyle
X\egalLoi U^{m}\bigg( \sum_{k=1}^{a+1}X^{(k)}+\sum_{k=a+2}^{S+1}Y^{(k)}\bigg)\\ \\
\displaystyle Y\egalLoi U^{m}\bigg(\sum_{k=1}^{c}X^{(k)}+\sum_{k=c+1}^{S+1}Y^{(k)} \bigg),
\end{array}
\right.
$$
where $U$ is uniform on $[0,1]$, $X^{(k)}$ and $Y^{(k)}$ are respective copies of $X$ and $Y$, all being independent of each other and of $U$.
Our System~\eqref{pointFixeDiscret} for the discrete time limit $W^{DT}$, slightly more complicated, is essentially of the same type.
These systems can be seen as natural generalizations of equations of type~\eqref{smoothingEq}, as set out in Neininger-R\"uschendorf~\cite{NeiRusaap}. 
%Adapting what has been done to study equations of type~\eqref{smoothingEq}, 
Section~\ref{sec:contractions} is devoted to the existence and the unicity of solutions of our systems by means of a contraction method (Theorems~\ref{th:solutionsSystemeDiscret} and~\ref{th:solutionsSystemeContinu}), leading to a characterization of $W^{DT}$ and $W^{CT}$ distributions.

Finally, in Section~\ref{sec:moments}, we take advantage of the fixed point systems again to give accurate bounds on the moments of $W^{CT}$ (Lemma~\ref{borne}).
Using this lemma, we show that the laws of $W^{DT}$ and $W^{CT}$ are determined by their moments (Corollary to Theorem~\ref{detMoment}).

%%%%%%%%%%%%%%%%%%%%%%%%%%%%%%%%
\section{Two-color P\'olya urn: definition and asymptotics}
\label{sec:model}
%%%%%%%%%%%%%%%%%%%%%%%%%

\subsection{Notations and asymptotics in discrete time}
\label{introDT}
%%%%%%%%%%%

Consider a two-color P\'olya-Eggenberger urn random process.
We adopt notations of the introduction:
the replacement matrix $R=\begin{pmatrix}a&b\\c&d\end{pmatrix}$ is assumed to have nonnegative entries, the integers $S$ as balance and $m$ as second smallest eigenvalue.
We assume $R$ to be non triangular, \emph{i.e.} that $bc\neq 0$;
this implies that $m\leq S-1$.
Moreover, the paper deals with \emph{large} urns which means that the ratio $\sigma =m/S$ is assumed to satisfy
$$
\sigma >\frac 12.
$$
We denote by $v_1$ and $v_2$ the vectors
\begin{equation}
\label{v}
v_1=\frac S{(b+c)}\begin{pmatrix}c\\b\end{pmatrix}
\hskip 1cm \hbox{ and } \hskip 1cm 
v_2 = \frac S{(b+c)}\begin{pmatrix}1\\-1\end{pmatrix};
\end{equation}
they are eigenvectors of the matrix $^t\! R$, respectively associated with the eigenvalues $S$ and~Ê$m$.
Let also $\left( u_1,u_2\right)$ be the dual basis
\begin{equation}
\label{u}
u_1(x,y) = \frac{1}{S}(x+y)
\hskip 1cm \hbox{ and } \hskip 1cm 
u_2(x,y)= \frac{1}{S}(bx-cy);
\end{equation}
$u_1$ and $u_2$  are eigenforms of $^t\! R$, respectively associated with the eigenvalues $S$ and~$m$.

\vskip 5pt
When the urn contains $\alpha$ white balls and $\beta$ black balls at (discrete) time $0$, the composition vector at time $n\in\g N$ is denoted by
$$
U_{(\alpha,\beta)}^{DT}(n).
$$
%It is known (see Janson \cite{Jan} and \cite{Pou08}) that the asymptotic behavior of the composition vector $U_{(\alpha,\beta)}^{DT}(n)$ when $n$ tends to $+\infty$ depends on the \emph{ratio} of the urn.
%When $\sigma \leq \frac 12$ the asymptotic behavior is Gaussian and one says that the urn is \emph{small}, though when $\sigma > \frac 12$, the urn is said \emph{large}.
Since the urn is assumed to be large, the asymptotics of its composition vector is given by the following result. %(see Janson~\cite{Jan}, or \cite{ChaPouSah,Pou08}).
\begin{Th}
\label{asympDT}
{\bf (Asymptotics of discrete time process,~\cite{Jan,Pou08})}

Let $\left( U_{(\alpha,\beta)}^{DT}(n)\right) _{n\in\g N}$ be a \emph{large} P\'olya urn discrete time process.
Then, when $n$ tends to infinity,
\begin{equation}
\label{asymptotiqueDiscrete}
U_{(\alpha,\beta)}^{DT}(n) = nv_1 + n^{\sigma} W_{(\alpha,\beta)}^{DT} v_2 + o( n^{\sigma} )
\end{equation}
where $v_1$ and $v_2$ are the non random vectors defined by \eqref{v}, $W_{(\alpha,\beta)}^{DT}$ is the real-valued random variable defined by 
\begin{equation}\label{martingalediscrete}
W_{(\alpha,\beta)}^{DT} := \lim_{n\rightarrow +\infty} \frac{1}{n^{\sigma}} u_2\left(U_{(\alpha,\beta)}^{DT}(n)\right)
\end{equation}
$u_2$ being defined in \eqref{u},  and where $o(\ )$ means almost surely and in any $L^p, p\geq 1$.
\end{Th}

\vskip 5pt
A proof of this result can be found in Janson~\cite{Jan} by means of embedding in continuous time, 
under an irreducibility assumption.
Another proof, which is valid in any case and that remains in discrete time is also given in~\cite{Pou08}.
The present paper is focused on the distribution of $W_{(\alpha,\beta)}^{DT}$ which appears in both proofs as the limit of a bounded martingale.
One remarkable fact that does not occur for small urns
(\emph{i.e.} when $\sigma\leq 1/2$) is that the distribution of $W_{(\alpha,\beta)}^{DT}$ actually depends on the initial composition vector $\left(\alpha,\beta\right)$.
For example, its expectations turns out to be
\be
\label{esperanceDT}
\g E\left( W_{(\alpha,\beta)}^{DT}\right)
=
\frac{\Gamma\left(\frac {\alpha +\beta}S\right)}{\Gamma\left(\frac {\alpha +\beta}S+\sigma\right)}\frac{b\alpha -c\beta}{S}.
\ee
This formula, explicitely stated in~\cite{ChaPouSah} can be shown by elementary means or using the convergent martingale
$$
\left(
\frac
{u_2\left( U_{(\alpha,\beta)}^{DT}(n)\right)}
{\displaystyle\prod _{0\leq k\leq n-1}\left( 1+\frac{\sigma}{k+\frac{\alpha +\beta}{S}}\right)}
\right) _{n\in\g N}.
$$
For more developments about this discrete martingale which is the essential tool in the discrete method for proving Theorem~\ref{asympDT}, see~\cite{Pou08}.

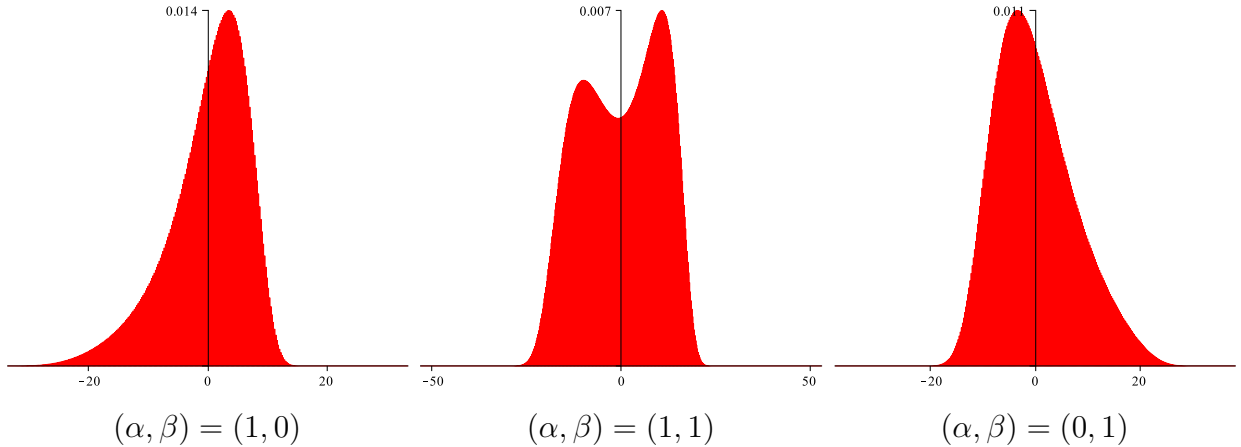
\begin{figure}
\begin{tikzpicture}
\draw (0,0) node{\includegraphics[width=55mm]{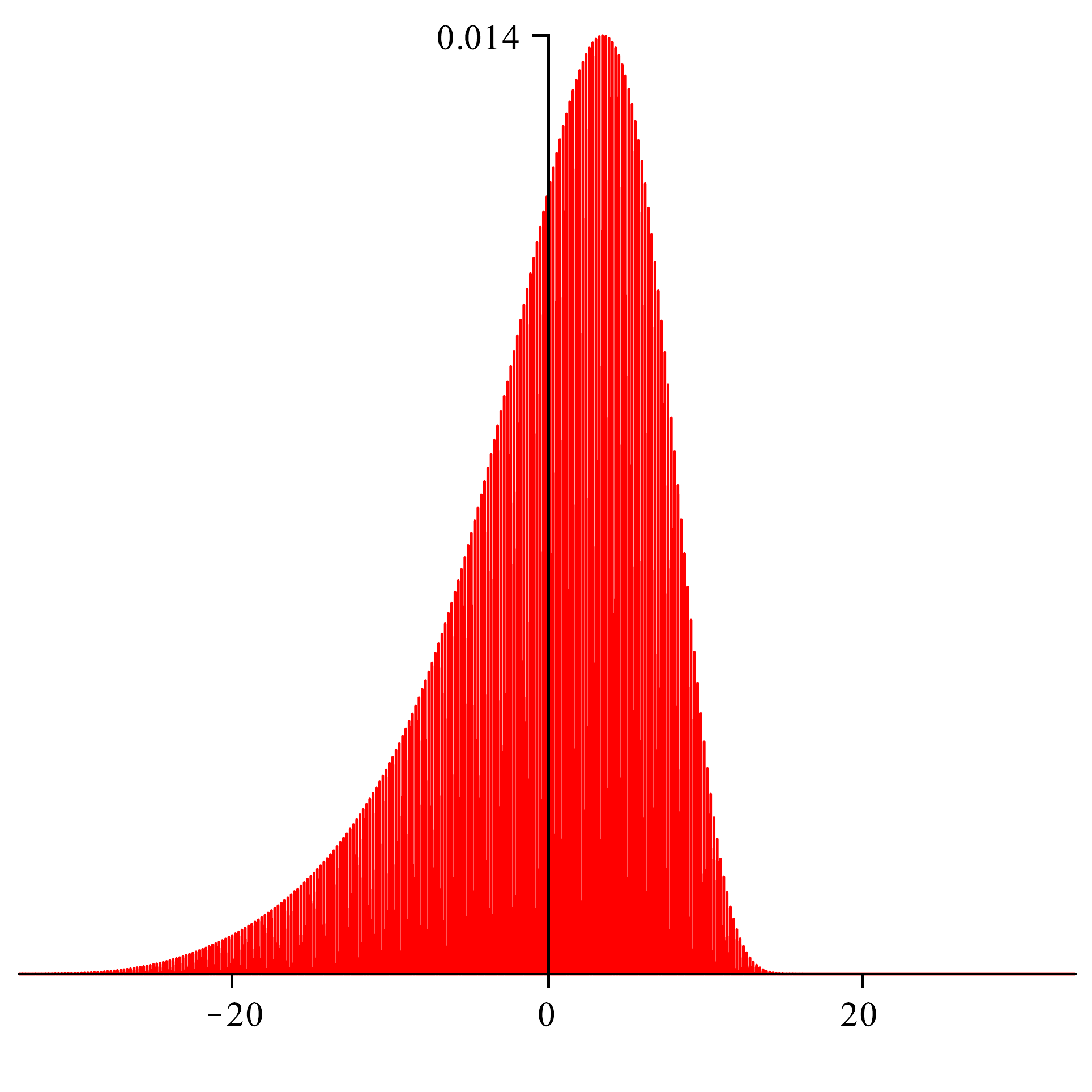}};
\draw (5.5,0) node{\includegraphics[width=55mm]{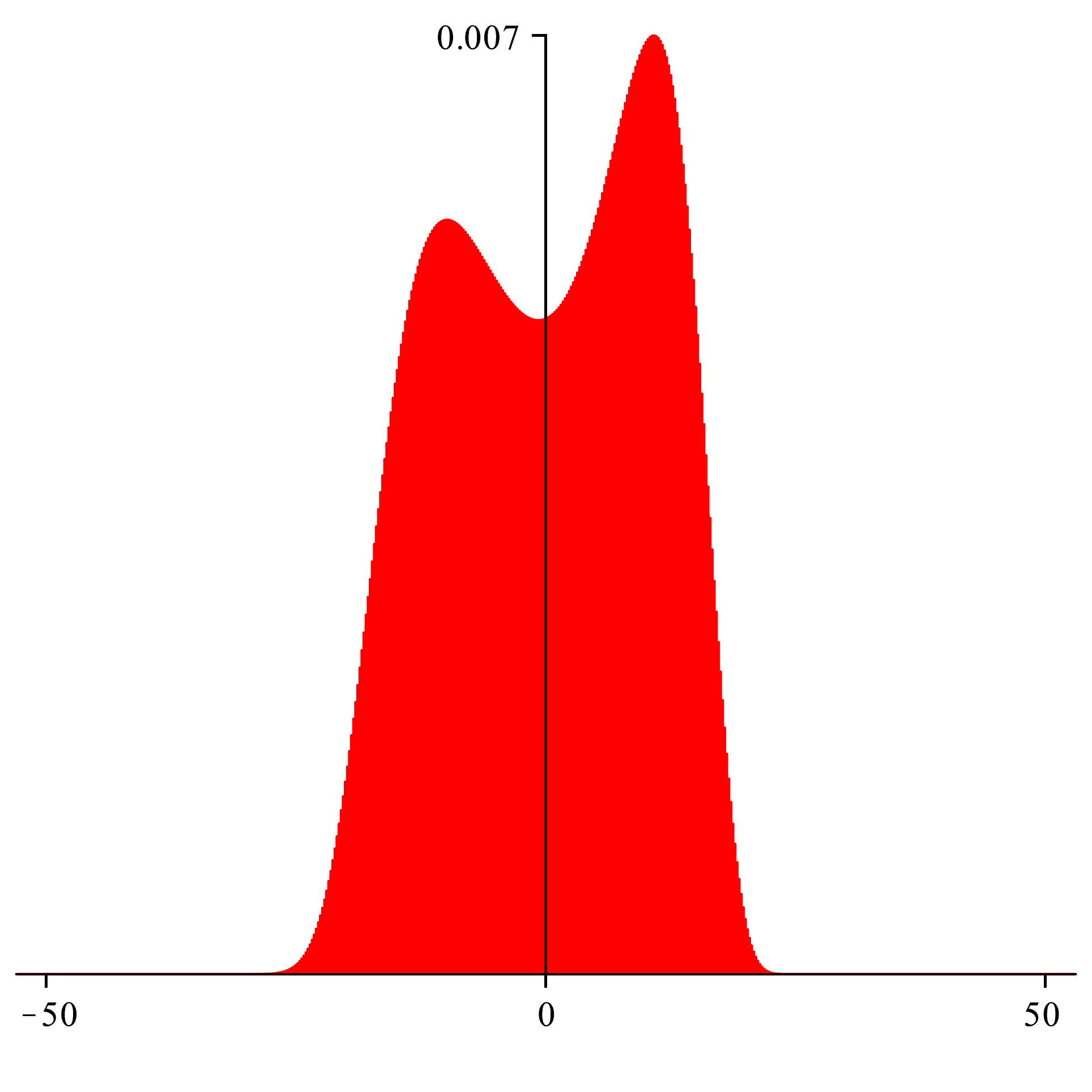}};
\draw (11,0) node{\includegraphics[width=55mm]{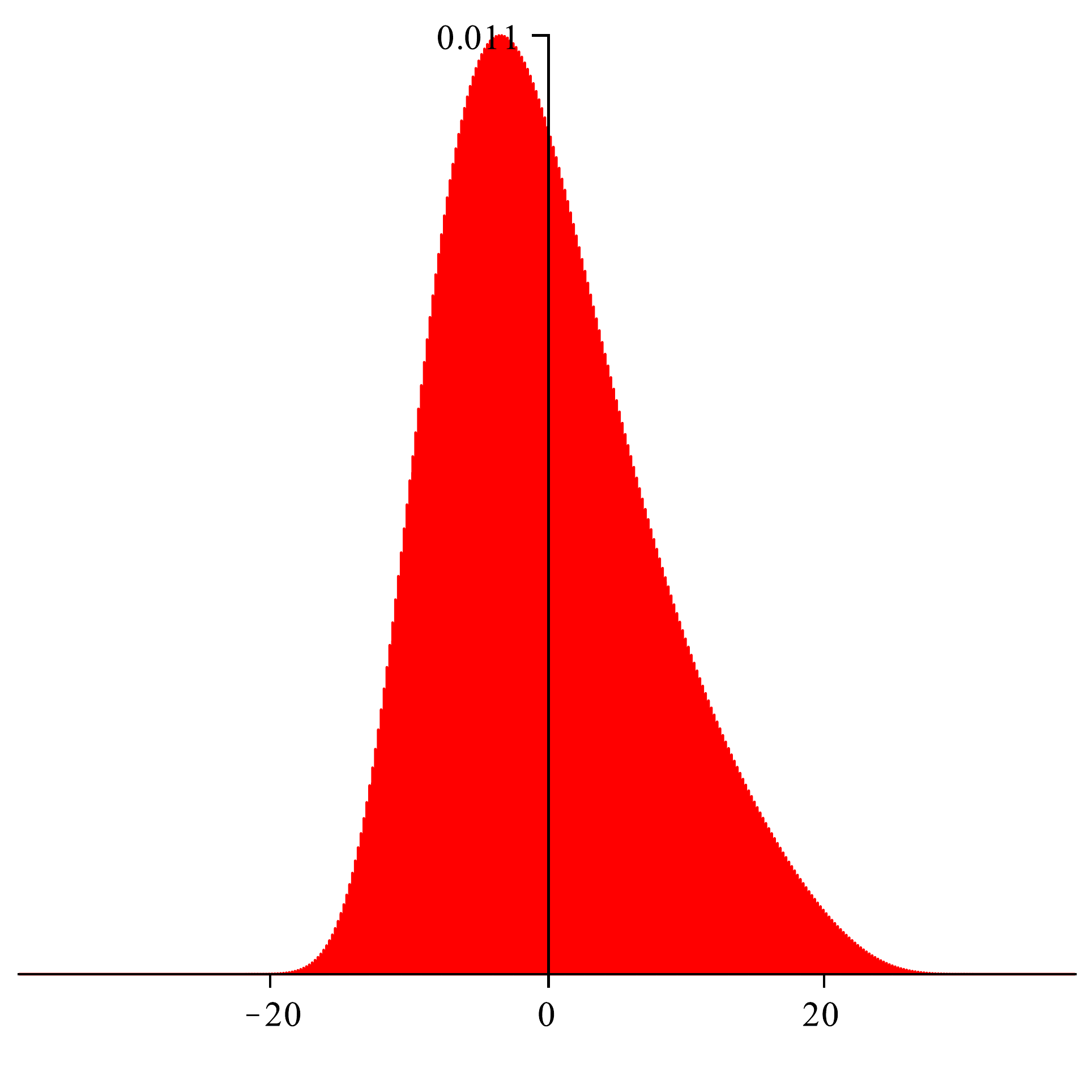}};
\draw (0,-3) node{$(\alpha ,\beta )=(1,0)$};
\draw (5.5,-3) node{$(\alpha ,\beta )=(1,1)$};
\draw (11,-3) node{$(\alpha ,\beta )=(0,1)$};
\end{tikzpicture}
\caption{\label{simul}starting from initial composition $(\alpha ,\beta )$, exact distribution of the number of red balls after $n=300$ drawings, centered around its mean and divided by $n^\sigma$.}
\end{figure}

\vskip 5pt
%The deterministic drift $nv_1$ in Formula~\eqref{asymptotiqueDiscrete} provides a good opportunity to approximate the distribution of $W^{DT}$.
%Indeed, it suffices to normalize the number of (say) red balls at time $n$ (first coordinate of the vector $U_{(\alpha,\beta)}^{DT}(n)$) by non random operations to get, when $n$ is large, an approximation of $W_{(\alpha,\beta)}^{DT}$.
The approach in analytic combinatorics makes easy to compute the probability generating function of the number of (say) red balls in the urn at finite time, by iteration of some suitable partial differential operator.
The treatment of P\'olya urns by analytic combinatorics is due to P. Flajolet and his co-authors and can be found in~\cite{FlaGabPek}.
Figure~\ref{simul} is the exact distribution of the (normalized) number of red balls after $300$ drawings, centered around its expectation.
The computations have been managed using Maple and concern the (large) urn with replacement matrix $R=\begin{pmatrix}18&2\\3&17\end{pmatrix}$ and respective initial compositions $(1,0)$, $(1,1)$ and $(0,1)$.

Some direct first observations can be made on these pictures.
For example, one gets an illustration of the decomposition formula~\eqref{decompDT} which states that the distribution of $U_{(1,1)}$ is decomposed as a weighted convolution of $U_{(1,0)}$'s and $U_{(0,1)}$'s.

%%%%%%%%%%%%%%%%%
\subsection{Embedding in continuous time}
%%%%%%%%%%%%%%%%%

Classically, the discrete time process is embedded in a continuous time multitype branching process; the idea of embedding discrete urn models in continuous time branching processes goes back
at least to Athreya and Karlin \cite{AK} and a description is given in Athreya and Ney \cite{AN}, Section 9.
The method has been revisited and developed by Janson \cite{Jan} and we summarize hereunder the results obtained in  \cite{ChaPouSah}.

We define the continuous time Markov  branching process
$$
\left( U_{(\alpha,\beta)}^{CT}(t)\right) _{t\in\g R_{\geq 0}}
$$
as being the embedded process of $\left( U_{(\alpha,\beta)}^{DT}(n)\right) _{n\in\g N}$.
It starts from the same initial condition
$U_{(\alpha,\beta)}^{CT}(0)=U_{(\alpha,\beta)}^{DT}(0)=(\alpha, \beta)$;
at any moment, each ball is equipped with an $\rond E xp(1)$-distributed\footnote{
For any positive real $a$, $\rond E xp(a)$ denotes the exponential distribution with parameter $a$.}  random clock, all the clocks being independent.
When the clock of a white ball rings, $a$ white balls and $b$ black balls are added in the urn;
when the ringing clock belongs to a black ball, one adds $c$ white balls and $d$ black balls, so that
the replacement rules are the same as in the discrete time urn process.

The important benefit of considering such a process comes from the \emph{independence} of the subtrees in the branching process. In the continuous-time urn process, each ball reproduces independently from the other balls. %In the discrete process indeed, the ``grapes'' of balls produced by replacement are not independent from each other.
%More details about the tree representation of the discrete-time process are given in in Section~\ref{sec:convolutionDT} 
The asymptotics of this process is given by the following theorem.%, proved in \cite{ChaPouSah} (a slightly weaker form being already stated in Janson~\cite{Jan}).

\begin{Th}
\label{continuousurn}
{\bf (Asymptotics of continuous time process, \cite{Jan,ChaPouSah})}

Let $\left( U_{(\alpha,\beta)}^{DT}(t)\right) _{t\geq 0}$ be a \emph{large} P\'olya urn continuous time process.
Then, when $t$ tends to infinity,
\begin{equation}
\label{asymptotiqueCT}
U_{(\alpha,\beta)}^{CT}\left( t\right)
=e^{St} \xi v_1\left( 1+o(1)\right)
+e^{mt}W_{(\alpha,\beta)}^{CT}v_2\left( 1+o(1)\right) ,
\end{equation}
where $v_1, v_2, u_1, u_2$ are defined in \eqref{v} and \eqref{u}, 
$\xi$ and $W_{(\alpha,\beta)}^{CT}$ are real-valued random variables defined by
$$
\xi =\lim _{t\to +\infty}e^{-St}u_1\left( U^{CT}(t)\right),
$$
$$
W_{(\alpha,\beta)}^{CT} =\lim _{t\to +\infty}e^{-mt}u_2\left( U_{(\alpha,\beta)}^{CT}(t)\right),
$$
all the convergences are almost sure and
in any ${\rm L}^p$-space, $p\geq 1$.
Furthermore, $\xi$ is {\rm Gamma}$\left(\frac{\alpha + \beta}{S}\right)$ distributed.%, where $u =\alpha +\beta$ is the total number of balls at time $0$.
\end{Th}

\vskip 5pt
Here again, the distribution of $W^{CT}$ depends on the initial composition of the urn.
For exemple, its expectation is
\be
\label{esperanceCT}
\g E\left( W_{(\alpha,\beta)}^{CT}\right)
=
\frac{b\alpha -c\beta}{S}.
\ee
as can be seen from the continuous-time martingale%~(\cite{Jan,ChaPouSah}).
$$
\Big(
e^{-mt}
u_2\left( U_{(\alpha,\beta)}^{DT}\left( t\right)\right)
\Big) _{t\geq 0}.
$$

Some properties of $W^{CT}$ are already known.
For example, it is supported by the whole real line $\g R$ and admits a density.
Moreover, this density is increasing on $\g R_{<0}$, decreasing on $\g R_{>0}$ and is not bounded in the neighbourhood of the origin.
Note that it is not an even function since $W^{CT}$ is not centered.
Finally, the characteristic function of $W^{CT}$ (\emph{i.e.} its Fourier transform) is infinitely differentiable but not analytic at the origin:
the domain of analyticity of $\g E\exp\left( zW^{CT}\right)$ is of the form
$\g C\setminus L_+\bigcup L_-$ where $L_+$ and $L_-$ are half-lines contained in~$\g R$, one of them being bordered at the the origin.
In particular, the exponential moment generating series of $W^{CT}$ has a radius of convergence equal to zero, due to a ramification \emph{and} a divergent series phenomenon as well.
All these properties are shown in~\cite{ChaPouSah}, based on the expression of this characteristic function in terms of the inverse of an abelian integral on the Fermat curve $x^m+y^m+z^m=0$.

%%%%%%%%%%%%%%%%%%%%%%%%%%%%%%%%%%%%
\subsection{Connexion discrete time/continuous time}
\label{sec:martingaleConnexion}

As in any embedding into continuous time of a Markov chain, the discrete time process and the continuous time one are connected by 
$$
\left(U^{CT}(\tau_n)\right)_{n\in\g N}= \left(U^{DT}\left( n\right)\right)_{n\in\g N}
$$
where 
$$
0=\tau_0 < \tau_1 < \dots  < \tau_n < \cdots
$$
are the jumping times of the continuous process. These random times are independent of the positions $U^{CT}(\tau_n)$.
The embedding for urn processes is widely studied in Janson~\cite{Jan}.
It is detailed in~\cite{ChaPouSah} in the special case of two-color P\'olya urns.
A dual formulation of this connexion is 
$$
\left(U^{CT}(t)\right)_{t\in\g R_{\geq 0}}= \left(U^{DT}\left( n(t)\right)\right)_{t\in\g R_{\geq 0}}
$$
where
$$
n(t):= \inf \{ n\geq 0, \tau_n \geq t \}
$$
is the number of drawings in the urn before time $t$.
After projection and normalization, these equalities provide the connexion between the limit variables $W_{(\alpha,\beta)}^{DT}$ and $W_{(\alpha,\beta)}^{CT}$:
\begin{equation}
\label{connection}
W_{(\alpha,\beta)}^{CT} \egalLoi \xi^{\sigma} \cdot W_{(\alpha,\beta)}^{DT}
%\hskip 1cm \hbox{ and } \hskip 1cm W_{(\alpha,\beta)}^{DT} \egalLoi \xi^{-\sigma} \cdot W_{(\alpha,\beta)}^{CT}
\end{equation}
where $\xi$ and $W_{(\alpha,\beta)}^{DT}$ are independent, $\xi$ being Gamma$\left(\frac{\alpha + \beta}{S}\right)$ distributed.
Note that the almost sure equality
$$
W_{(\alpha,\beta)}^{DT} =\xi^{-\sigma} \cdot W_{(\alpha,\beta)}^{CT}
$$
holds as well, the variables $\xi$ and $W_{(\alpha,\beta)}^{CT}$ being however \emph{not} independent.

%%%%%%%%%%%%%%%%%%%%%%%%%%%%%%%%%%%%
%%%%%%%%%%%%%%%%%%%%%%%%%%%%%%%%%%%%
\section{Decomposition properties}
\label{sec:decomposition}

This section emphasizes the underlying tree structure of the urn process. This obvious vision is indeed the key in the following two decompositions: first, we reduce the study of $W_{(\alpha ,\beta )}$ to the study of $W_{(1,0)}$ and
$W_{(0,1)}$, called later on $X$ and $Y$ respectively, to lighten the notations. Second, in Section \ref{fixedpointDT}, we exploit a ``divide-and-conquer'' property to deduce a system of fixed point equations on $X$ and $Y$. The reasoning is detailed in discrete time. It is much more straightforward in continuous time, since the decomposition is contained inside the branching property.
Detailed in \cite{ChaPouSah}, the continuous case is  briefly recalled in Section~\ref{convolutionCT}.

The natural question ``is it possible to deduce the DT-system from the CT-system and conversely'' is partially adressed in Section \ref{sec:connexion-systems}.

%%%%%%%%%%%%%%%%%%%%%%%%%%%%%%%%%%%%%
\subsection{Tree structure in discrete time}
\label{sec:convolutionDT}
%%%%%%%%%%%%%%%%%%%%%%%%%%%%%%%%%%%%%%

In this section dealing with the \emph{discrete time} process, we skip the index $DT$ when no confusion is possible. 

Let us make precise the tree structure of the urn process: a forest $(\rond T_n)$ grows at each drawing from the urn. At time $0$ the forest is reduced to $\alpha$ red nodes and $\beta$ black nodes, which are the roots of the forest trees. At time $n$, each leaf in the forest represents a ball in the urn.
When a leaf is chosen (a ball is drawn), it becomes an internal node and gives birth to $(a+1)$ red leaves and $b$ black leaves, or $c$ red leaves and $(d+1)$ black leaves, according to the color of the chosen leaf.

%- we call ``children'' of a ball $u$ the balls which are created when drawing the ball $u$ and 

%- we call ``subtree'' of a ball $u$ the set of  $u$ and its ``descendants''.

The dynamics of the urn process was described saying ``at each time $n$, a ball is uniformly chosen in the urn''. It becomes ``a leaf is uniformly chosen among the leaves of the forest''. This forest therefore appears as a non binary colored generalization of a binary search tree. 

For example, take the following urn with 
$R=\left(
\begin{array}{cc}
6&1\\
2&5
\end{array}
\right)$ as replacement matrix (it is a large urn) and start from $\alpha = 3$ red balls and  $\beta = 2$ black balls.  Below is a possible configuration after $3$ drawings.

\vskip 1cm

\begin{center}
\begin{tikzpicture}
[level 1/.style={sibling distance=6mm},
level 2/.style={sibling distance=4mm},
level 3/.style={sibling distance=5mm}]
\tikzstyle{noeud}=[circle,fill=black]
\tikzstyle{noire}=[circle,fill=black]
\tikzstyle{rouge}=[circle,fill=red]
\tikzstyle{rougevide}=[circle,draw=red,line width=2pt]
\tikzstyle{noirevide}=[circle,draw=black,line width=2pt]
\node [rouge] at (0,0) {};
\node [rouge] at (3,0) {};
\node [rougevide] at (6,0) {}
	child { node[rouge] {}}
	child { node[rouge] {}}
	child { node[rouge] {}}
	child { node[rouge] {}}
	child { node[rouge] {}}
	child { node[rouge] {}}
	child { node[rouge] {}}
	child { node[noire] {}}
	;
\node [noire] at (9,0) {};
\node [noirevide] at (12,0) {}
	child { node[rouge] {}}
	child { node[rouge] {}}
	child { node[noirevide] {}
		child { node[rouge] {}}
		child { node[rouge] {}}
		child { node[noire] {}}
		child { node[noire] {}}
		child { node[noire] {}}
		child { node[noire] {}}
		child { node[noire] {}}
		child { node[noire] {}}
	}
	child { node[noire] {}}
	child { node[noire] {}}
	child { node[noire] {}}
	child { node[noire] {}}
	child { node[noire] {}}
	;
\end{tikzpicture}
\end{center}

\medskip

Initial red balls are numbered from $1$ to $\alpha$ and initial black balls from $(\alpha + 1)$ to $(\alpha + \beta)$. The following figure represents the forest coming from these initial balls.

\begin{center}
\begin{tikzpicture}[scale=1]
\def\rayon{0.3}
\def\hauteur{3}
\def\blanche(#1,#2){\draw (#1,#2) circle (\rayon);}
\def\noire(#1,#2){\draw [line width=2pt](#1,#2) circle (\rayon);}
\def\rouge(#1,#2){\draw [red, line width=2pt] (#1,#2) circle (\rayon);}
\def\sousarbre(#1,#2){\fill [color=gray!50] (#1,#2-\rayon)--++(-2*\rayon,-\hauteur)--++(4*\rayon,0) ;}%sous une boule
\rouge(0,10)
\sousarbre(0,10)
\rouge(5*\rayon,10)
\sousarbre(5*\rayon,10)
\draw [dotted,line width=2pt] (8*\rayon ,10)--(12*\rayon ,10) ;
\rouge(15*\rayon,10)
\sousarbre(15*\rayon,10)
\noire(20*\rayon,10)
\sousarbre(20*\rayon,10)
\draw [dotted,line width=2pt] (23*\rayon ,10)--(27*\rayon ,10) ;
\noire(30*\rayon,10)
\sousarbre(30*\rayon,10)

\draw (7.5*\rayon,11) node{$\alpha$} ;
\draw (0,10.6)--(0,11)--(6.5*\rayon,11) ;
\draw (8.5*\rayon,11)--(15*\rayon,11)--(15*\rayon,10.6) ;
\draw (25*\rayon,11) node{$\beta$} ;
\draw (20*\rayon,10.6)--(20*\rayon,11)--(24*\rayon,11) ;
\draw (26*\rayon,11)--(30*\rayon,11)--(30*\rayon,10.6) ;
\end{tikzpicture}
\end{center}

\medskip
For any $n\geq 0$ and $k\in\{ 1,\dots ,\alpha +\beta\}$, denote by $D_k(n)$ the number of leaves of the $k$-th tree in the forest at time $n$.
Thus, at time~$n$, the number of drawings in the $k$-th tree is $\frac{D_k(n)-1}{S}$.
This numbers represents the \emph{time} inside this $k$-th tree.

Remember that the balls of the whole urn are uniformly drawn at any time and notice that at each drawing in the $k$-th tree, $D_k(n)$ increases by $S$:
the random vector $D(n)=\left( D_1(n),\dots ,D_{\alpha +\beta}(n)\right)$ has exactly the same distribution as the composition vector at time $n$ of an $(\alpha +\beta)$-color P\'olya urn process having $SI_{\alpha +\beta}$ as replacement matrix and $(1,\dots ,1)$ as initial composition vector.

Gathering these arguments, the distribution of $U_{(\alpha ,\beta )}(n)$ can be described the following way:
consider simultaneously

{\it (i)} an original $(\alpha +\beta)$-color urn process
%$D(n)=\left( D_1(n),\dots ,D_{\alpha +\beta}(n)\right)$
$D=\left( D_1,\dots ,D_{\alpha +\beta}\right)$
having $SI_{\alpha +\beta}$ as matrix replacement and $(1,\dots ,1)$ as initial condition;

{\it (ii)} for any $k\in\{ 1,\dots ,\alpha\}$, an urn process
%$\left( U_{(1,0)}^{(k)}(n)\right) _n$
$U_{(1,0)}^{(k)}$
having $R$ as replacement matrix and $(1,0)$ as initial condition;

{\it (iii)} for any $k\in\{ \alpha +1,\dots ,\alpha +\beta\}$, an urn process
%$\left( U_{(0,1)}^{(k)}(n)\right) _n$
$U_{(0,1)}^{(k)}$
having $R$ as replacement matrix and $(0,1)$ as initial condition,

all these processes being independent of each other.
Then, the process
$U_{(\alpha ,\beta )}=\left( U_{(\alpha ,\beta )}(n)\right) _n$
has the same distribution as the process defined by the sum of the $U_{(1,0)}^{(k)}$ and of the $U_{(0,1)}^{(k)}$ at respective times $\frac{ D_k(n)-1}{S}$.
%Gathering these arguments, it appears that $U_{(\alpha ,\beta )}(n)$ has the same distribution as the sum of $\alpha$ processes $U_{(1,0)}^{(k)}$ $(1\leq k\leq\alpha )$ and $\beta$ processes $U_{(0,1)}^{(k)}$ $(\alpha +1\leq k\leq\alpha +\beta )$ independent processes each being taken at time $\frac{ D_k(n)-1}{S}$, where $\left( D_1(n),\dots ,D_{\alpha +\beta}(n)\right)$ is an original P\'olya urn process:
In other words, for any $n\geq 0$,
\be
\label{convolutionFiniteTime}
U_{(\alpha ,\beta )}\left( n\right)\egalLoi
\sum _{k=1}^\alpha U_{(1,0)}^{(k)}\big( \hbox{$\scriptsize \frac{ D_k(n)-1}{S}$}\big)
+\sum _{k=\alpha +1}^{\alpha +\beta}U_{(0,1)}^{(k)}\big( \hbox{$\scriptsize \frac{ D_k(n)-1}{S}$}\big)
\ee
where the $U_{(1,0)}^{(k)}$ and the $U_{(0,1)}^{(k)}$ are respective copies of the random vector
processes $U_{(1,0)}$ and $U_{(0,1)}$, all being independent of each other and of $D$.

\medskip
%Consequently, the
The following claim is a direct consequence of Proposition~\ref{proDirichlet} in Section~\ref{sec:appendix}.

{\bf Claim\ }\hskip 10pt
{\it
When $n$ goes off to infinity, $\frac 1{nS}\left( D_1(n),\dots ,D_{\alpha +\beta}(n)\right)$ converges
almost surely to a $Dirichlet\left( \frac 1S,\dots ,\frac 1S\right)$- distributed random vector, denoted by
$Z=(Z_1,\dots ,Z_{\alpha +\beta})$.
}

\medskip

Notice that for any $k$, $D_k(n)$ tends almost surely to $+\infty$ when $n$ tends to infinity.
Starting from Equation~(\ref{convolutionFiniteTime}), dividing by $n^\sigma$,
taking the image by the second projection $u_2$ (notations of Section~\ref{introDT}) and passing to the (almost sure) limit $n\to\infty$ thanks to Theorem~\ref{asympDT}, one obtains the following theorem.
\begin{Th}
\label{thConvolution}
For any $(\alpha ,\beta )\in\g N^2\setminus (0,0)$, let $W_{(\alpha ,\beta )}$ be the limit distribution
of a large two-color \emph{discrete} time P\'olya urn process with ratio $\sigma$ and initial condition
$(\alpha ,\beta )$.
Then,
\be
\label{decompDT}
W_{(\alpha ,\beta )}
\egalLoi
\sum _{k=1}^{\alpha}
Z_k^\sigma W_{(1,0)}^{(k)}+\sum _{k=\alpha +1}^{\alpha +\beta}Z_k^\sigma W_{(0,1)}^{(k)}
\ee
where

(i) $Z=(Z_1,\dots ,Z_{\alpha +\beta})$ is a Dirichlet distributed random vector, with parameters
$(\frac 1S,\dots ,\frac 1S)$;

(ii) the $W_{(1,0)}^{(k)}$ and the $W_{(0,1)}^{(k)}$ are respective copies of $W_{(1,0)}$ and
$W_{(0,1)}$, all being independent of each other and of $Z$.
\end{Th}

Notice that any $Z_k$ is $Beta (\frac 1S,\frac{\alpha +\beta -1}{S})$-distributed (see Section~\ref{sec:appendix}).

%%%%%%%%%%%%%%%%%%%%%%%%%%%%%%%
\subsection{Discrete time fixed point equation}
\label{fixedpointDT}
%%%%%%%%%%%%%%%%%%%%%%%%%%%%%%%%

Theorem~\ref{thConvolution} shows that the limit distribution of a large urn process starting with any
initial composition can be written as a function of two ``elementary'' particular laws, namely the laws  of $W^{DT}_{(1,0)}$
and $W^{DT}_{(0,1)}$.
The present section gives a characterisation of these two distributions by means of a fixed point
equation.

Let $(U(n))_{n\geq 0}$ be a two-color P\'olya urn process, with all the notations of
Section~\ref{introDT}.
In order to simplify the notations, denote 
\be
\label{defXY}
\left\{
\begin{array}{l}
\displaystyle X:=
W^{DT}_{(1,0)}=\lim_{n\rightarrow +\infty} u_2 \left( \frac{U_{(1,0)}(n)}{n^{\sigma}}\right) \\ \\
\displaystyle Y:=
W^{DT}_{(0,1)}=\lim_{n\rightarrow +\infty} u_2 \left( \frac{U_{(0,1)}(n)}{n^{\sigma}}\right)
\end{array}
\right.
\ee
Focus now on the study of $U_{(1,0)}(n)$. At time $1$ the composition of the urn is deterministic: there are $(a+1)$ red balls and $b$ black balls. Exactly like in Section \ref{sec:convolutionDT}, the tree structure of the urn appears, with a forest starting from $(a+1)$ red balls and $b$ black balls. In the same example with replacement matrix 
$R=\left(
\begin{array}{cc}
6&1\\
2&5
\end{array}
\right)$,
this fact is illustrated by the following figure:
\begin{center}
\tikzstyle{bouleblanche}=[circle,draw]
\tikzstyle{boulenoire}=[circle,draw=black,line width=2pt]
\tikzstyle{boulerouge}=[circle,draw=red,line width=2pt]

\begin{tikzpicture}
[level 1/.style={sibling distance=19mm}]
\def\rayon{0.35}
\def\hauteur{4}
\def\sousarbre(#1,#2){\fill [color=gray!50] (#1,#2-\rayon)--++(-2*\rayon,-\hauteur)--++(4*\rayon,0) ;}%sous une boule
\node [boulerouge] {}
	child {node [boulerouge] {} }
	child {node [boulerouge] {} }
	child {node [boulerouge] {} }
	child {node [boulerouge] {} }
	child {node [boulerouge] {} }
	child {node [boulerouge] {} }
	child {node [boulerouge] {} }
	child {node [boulenoire] {} };
\foreach \n in {-3,-2,-1,0,1,2,3,4} \sousarbre(-0.95+1.9*\n,-1.4);
\end{tikzpicture}
\end{center}

\medskip

For any $n\geq 1$, denote by $J_k(n)$ the number of leaves at time $n$ of the $k$-th subtree.
Then, at time $n$, the number of drawings in the $k$-th subtree is $\frac{J_k(n)-1}{S}$ so that, as in Section \ref{sec:convolutionDT}, one gets the equation in distribution
\be
\label{divideConquer}
U_{(1,0)}\left( n\right)\egalLoi
\sum_{k=1}^{a+1} U_{(1,0)}^{(k)} \big( \hbox{$\scriptsize \frac{ J_k(n)-1}{S}$}\big)
+\sum_{k=a+2}^{S+1} U_{(0,1)}^{(k)} \big( \hbox{$\scriptsize \frac{ J_k(n)-1}{S}$}\big)
\ee
where the $U_{(1,0)}^{(k)}$ and the $U_{(0,1)}^{(k)}$ are respective copies of the random vector
processes $U_{(1,0)}$ and $U_{(0,1)}$, all being independent of each other and of the $J_k$'s.
Besides, the random vector $\left( J_1(n),\dots ,J_{S+1}(n)\right)$ is exactly distributed like the composition vector at time $(n-1)$ of an $(S+1)$-color P\'olya urn process having $SI_{S+1}$ as replacement matrix
and $(1,\dots ,1)$ as initial composition vector, so that, by Proposition~\ref{proDirichlet} in Section~\ref{sec:appendix},
$$
\frac 1{nS}\Big( J_1(n),\dots ,J_{S+1}(n)\Big) \conv
V=\left(V_1,\dots ,V_{S+1}\right)
$$
almost surely, the random vector $V$ being $Dirichlet\left(\frac 1S,\dots ,\frac 1S\right)$-distributed.
Like in Section~\ref{sec:convolutionDT}, divide Equation~(\ref{divideConquer}) by $n^\sigma$, take the image by the second projection~$u_2$ and pass to the limit $n\to\infty$ using Theorem~\ref{asympDT}.
This leads to the following theorem.

\begin{Th}
\label{th:systemDT}
As defined just above by~(\ref{defXY}), let $X$ and $Y$ be the elementary limit laws of a large two-color \emph{discrete} time P\'olya urn process with replacement matrix
$\left(
\begin{array}{cc}
a&b\\c&d
\end{array}
\right)$,
balance $S=a+b=c+d$ and ratio~$\sigma >\frac 12$.
Then, $X$ and $Y$ satisfy the distributional equations system
\be
\label{pointFixeDiscret}
%\label{pointFixeDiscret}
\left\{
\begin{array}{l}
\displaystyle
X\egalLoi
\sum_{k=1}^{a+1} V_k ^{\sigma} X^{(k)}
+\sum_{k=a+2}^{S+1}  V_k^{\sigma} Y^{(k)}
\\ \\
\displaystyle
Y\egalLoi
\sum_{k=1}^{c} V_k^{\sigma} X^{(k)}
+\sum_{k=c+1}^{S+1} V_k^{\sigma} Y^{(k)}
\end{array}
\right.
\ee
where

(i) $V=(V_1,\dots ,V_{S+1})$ is a Dirichlet distributed random vector, with parameters
$(\frac 1S,\dots ,\frac 1S)$;

(ii) the $X^{(k)}$ and the $Y^{(k)}$ are respective copies of $X$ and $Y$, all being independent of
each other and of $V$.
\end{Th}

Notice that any $V_k$ is distributed like a random variable $U^S$, $U$ being uniformly distributed on $[0,1]$.
Equivalently, $V_k^{\sigma}$ is distributed like $U^m$ (notations of Section~\ref{introDT}).

%%%%%%%%%%%%%%%%%%%%%%%%%%%%%%%%%%%%%
\subsection{Decomposition properties in continuous time}
\label{convolutionCT}
%%%%%%%%%%%%%%%%%%%%%%%%%%%%%%%%%%%%%%

Remember that $\left( U^{CT}(t)\right) _t$ is a continuous time branching process. 
Thanks to the branching property, the decomposition properties of this process are somehow automatic. First, 
$$
U_{(\alpha,\beta)}^{CT}\left( t\right) =
\left[\alpha \right] U_{(1,0)}^{CT}\left( t\right) +\left[ \beta\right] U_{(0,1)}^{CT}\left( t\right) ,
$$
where the notation $[n]X$ means the sum of $n$ independant random variables having the same distribution as $X$. Consequently, passing to the limit when $t\rightarrow +\infty$ after normalization and projection yields
\begin{equation}
\label{decompCT}
W^{CT}_{(\alpha,\beta)}=
\left[ \alpha\right] W^{CT}_{(1,0)}+\left[ \beta\right] W^{CT}_{(0,1)}.
\end{equation}
This convolution formula expresses how the limit law $W^{CT}$ is decomposed in terms of elementary limit laws $W^{CT}_{(1,0)}$ and $W^{CT}_{(0,1)}$.
It corresponds to the discrete time decomposition shown in Theorem~\ref{thConvolution}.

\vskip 5pt
Now start from one red ball or from one black ball, and apply again the branching property at the first splitting time.
As before, define $X^{CT}$ and $Y^{CT}$ by
\be
\label{defXYcontinuous}
\left\{
\begin{array}{l}
\displaystyle X^{CT}:= W^{CT}_{(1,0)} = \lim _{t\to +\infty}e^{-mt}u_2\left( U_{(1,0)}^{CT}(t)\right), \\ \\
\displaystyle Y^{CT}:= W^{CT}_{(0,1)} = \lim _{t\to +\infty}e^{-mt}u_2\left( U_{(0,1)}^{CT}(t)\right) .
\end{array}
\right.
\ee
Then, with the above Theorem \ref{continuousurn}, one gets the following result.%, already mentioned in Janson~\cite{Jan} and in~\cite{ChaPouSah}.

\begin{Th}[\cite{Jan, ChaPouSah}]
\label{th:systemCT}
Let $X=X^{CT}$ and $Y=Y^{CT}$ be the elementary limit laws of a large two-color \emph{continuous} time P\'olya urn process with replacement matrix
$\left(
\begin{array}{cc}
a&b\\c&d
\end{array}
\right)$,
balance $S=a+b=c+d$ and ratio~$\sigma >\frac 12$, as defined just above by~(\ref{defXYcontinuous}).
Then, $X$ and $Y$ satisfy the distributional equations system
\begin{equation}
%\label{pointFixeContinu}
\label{pointFixeContinu}
\left\{
\begin{array}{l}
\displaystyle
X\egalLoi U^{m}\bigg( \sum_{k=1}^{a+1}X^{(k)}+\sum_{k=a+2}^{S+1}Y^{(k)}\bigg)\\ \\
\displaystyle Y\egalLoi U^{m}\bigg(\sum_{k=1}^{c}X^{(k)}+\sum_{k=c+1}^{S+1}Y^{(k)} \bigg),
\end{array}
\right.
\end{equation}
where $U$ is uniform on $[0,1]$,
where $X$,  $X^{(k)}$ and $Y$,  $Y^{(k)}$ are respective copies of $X^{CT}$ and $Y^{CT}$, all being independent of
each other and of $U$.
\end{Th}

\begin{Rem}
As mentioned above, it is shown in~\cite{ChaPouSah} that $X^{CT}$ (and $Y^{CT}$) admit densities.
The proof is based on the computation of the Fourier transform of $X^{CT}$ in terms of the inverse of an abelian integral on a Fermat curve.
This method is specific to $2$-color urn processes.
Theorems~\ref{th:systemDT} and~\ref{th:systemCT} give a new way of proving this fact by means of techniques that can be adapted from Liu's method (see~\cite{Liu99} for example).
This alternative method provides a perspective (adressed in a forthcoming paper):
it can be applied to show that the limit laws of $d$-color large urns admit densities as well.
\end{Rem}

%%%%%%%%%%%%%%%%%%%%%%%%
\subsection{Connexion between continuous-time and discrete-time systems}
\label{sec:connexion-systems}

In Section~\ref{sec:martingaleConnexion}, we described the connexion between the limit laws of large urns in discrete and continuous time, called the \emph{martingale connexion}.
It was seen as a consequence of the embedding into continuous time of the initial discrete time Markov chain defining the urn process.
In this paragraph, we show how one can deduce the solutions of the continuous time System~\eqref{pointFixeContinu} from the solutions of the discrete time System~\eqref{pointFixeDiscret}.

\begin{Prop}
\label{connexionSystems}
Let $X$ and $Y$ be solutions of the distributional System~\eqref{pointFixeDiscret} and let $\xi$ be a $Gamma$-distributed random variable with parameter $\frac 1S$, independent of $X$ and $Y$.
Then, $\xi ^\sigma X$ and $\xi ^\sigma Y$ are solutions of the distributional System~\eqref{pointFixeContinu}.
\end{Prop}

The assertion of Proposition~\ref{connexionSystems} is a consequence of the following lemma which is an elementary result in probability theory.

\begin{Lem}
Consider the two following distributional equations with unknown real-valued random variables $X, X_1,\dots ,X_{S+1}$.

1- \emph{Equation D:}
$$
X\egalLoi\sum _{1\leq k\leq S+1}V_k^\sigma X_k
$$
where ${\mathbf V}=\left( V_1,\dots ,V_{S+1}\right)$ is a Dirichlet-distributed random vector with parameter $\left(\frac 1S,\dots ,\frac 1S\right)$.

\vskip 5pt
2- \emph{Equation C:}
$$
X\egalLoi V^\sigma\sum _{1\leq k\leq S+1}X_k
$$
where $V$ is a $Beta$-distributed random variable with parameter $\left(\frac 1S,1\right)$ (in other words, $V^{1/S}$ is uniformly distributed on $[0,1]$).
\vskip 3pt
Let $V, \xi _1,\dots ,\xi _{S+1}$ be independent random variables, the $\xi _k$'s being
$Gamma\left(\frac 1S\right)$-distributed and $V$ being $Beta\left(\frac 1S,1\right)$-distributed.
Denote
$$
\displaystyle\xi :=V\sum _{1\leq j\leq S+1}\xi _j
$$
and, for any $k\in\{ 1,\dots S+1\}$,
$$
V_k:=\frac{\xi _k}{\displaystyle\sum _{1\leq j\leq S+1}\xi _j}.
$$

Then,
\begin{enumerate}[(i)]
\item
the random variable $\xi$ is $Gamma\left(\frac 1S\right)$-distributed;

\item
the random vector $\left( V_1,\dots ,V_{S+1}\right)$ is independent of $\xi$ and Dirichlet-distributed with parameter $\left(\frac 1S,\dots ,\frac 1S\right)$.

\item
if $X, X_1,\dots ,X_{S+1}$ satisfy Equation D, then $\xi ^\sigma X,\xi _1^\sigma X_1,\dots ,\xi _{S+1}^\sigma X_{S+1}$ satisfy Equation~C.
\end{enumerate}
\end{Lem}

\pff
{\it (i)}
This can be seen for example by computation of moments ($Beta$ and $Gamma$ distributions are moment-determined):
the $p$-th moment of a $Gamma (\alpha )$ distribution is
$\frac{\Gamma (\alpha +p)}{\Gamma (\alpha )}$
and the $p$-th moment of a $Beta (\alpha ,\beta)$ distribution is
$\frac{\Gamma (\alpha +p)\Gamma (\beta )}{\Gamma (\alpha +\beta +p)}$
where $\Gamma$ is Euler Gamma function.
Moreover, the sum of independent $Gamma(\alpha _1),\dots,Gamma(\alpha _d)$-distributed random variables is $Gamma(\alpha _1+\dots +\alpha _d)$-distributed.
Assertion {\it (i)} is a direct consequence of these facts.

{\it (ii)}
Classically, $\left( V_1,\dots ,V_{S+1}\right)$ is Dirichlet distributed and independent of the sum
$\sum _{1\leq j\leq S+1}\xi _j$.
For a proof of this result, see for example Chaumont and Yor~\cite{ChaumontYor}.
Since $V$ is independent of the $\xi _k$'s, the random variable $\xi$ is also independent of $\left( V_1,\dots ,V_{S+1}\right)$.

{\it (iii)}
Suppose that $X, X_1,\dots ,X_{S+1}$ satisfy Equation D.
Then $X\egalLoi\sum _{1\leq k\leq S+1}V_k^\sigma X_k$.
Multiplying the equality by the random variable $\xi ^\sigma$ leads to the distributive relation
$$
\xi ^\sigma X
\egalLoi
\left(\frac{\xi}{\displaystyle\sum _{1\leq j\leq S+1}\xi _j}\right) ^\sigma\sum _{1\leq k\leq S+1}\xi _k^\sigma X_k
=V^\sigma\sum _{1\leq k\leq S+1}\xi _k^\sigma X_k
$$
which makes the proof complete.
\QED

%%%%%%%%%%%%%%%%%%%%%%%%%%%
\subsection{Densities}
\label{sec:shape}

As shown in~\cite{ChaPouSah}, the law of $W_{(\alpha,\beta)}^{CT}$ turns out to be absolutely continuous with regard to Lebesgue measure on~$\g R$.
In this section, the same property is deduced for $W_{(1,0)}^{DT}$ and $W_{(0,1)}^{DT}$ from the fixed point Equation~\eqref{pointFixeDiscret}.
The decomposition property~\eqref{decompDT} then implies that any $W_{(\alpha,\beta)}^{DT}$ also admits a density.
The observations made in Section \ref{introDT} on Figure~\ref{simul} can be seen as a first approximation of the shape of the density of $W_{(\alpha,\beta)}^{DT}$.

The method we use is widely inspired from Q. Liu papers~\cite{Liu99} and~\cite{Liu01}.
See also ~\cite{ChaLiuPouContinu,ChaLiuPouDiscret} for an argumentation of the same vein for complex-valued probability measures.
Applied to fixed point Equation~\eqref{pointFixeContinu}, this method provides a second proof for the absolute continuity of $W_{(\alpha,\beta)}^{CT}$;
details are left to the reader.

For the whole Section~\ref{sec:shape}, we denote
$$
X:=W_{(1,0)}^{DT}
{\rm ~and~~}
Y:=W_{(0,1)}^{DT}.
$$
Let $\varphi _X$ and $\varphi _Y$ be the Fourier transforms of $X$ and $Y$:
for any $t\in\g R$,
$$
\varphi _X(t)=\g E\left( e^{itX}\right)
{\rm ~and~~}
\varphi _Y(t)=\g E\left( e^{itY}\right).
$$

\begin{Th}
\label{densiteXYDT}
As defined just above, let $X$ and $Y$ be the elementary limit distributions of a large two-color \emph{discrete} time P\'olya urn process with replacement matrix
$R=\left(
\begin{array}{cc}
a&b\\c&d
\end{array}
\right)$.
As in the whole paper, let $S=a+b=c+d$ and $m=a-c=d-b\in ]\frac S2,S[$ both eigenvalues of~$R$.
Then,

\begin{enumerate}[(i)]
\item
the support of $X$ and $Y$ is the whole real line $\g R$;

\item
for any $\rho\in ]0,\frac{a+1}m[$, there exists $C>0$ such that for any $t\in\g R\setminus\{ 0\}$,
$$
\left|\varphi _X(t)\right|\leq\frac{C}{\left| t\right| ^\rho};
$$

\item
for any $\rho\in ]0,\frac{d+1}m[$, there exists $C>0$ such that for any $t\in\g R\setminus\{ 0\}$,
$$
\left|\varphi _Y(t)\right|\leq\frac{C}{\left| t\right| ^\rho};
$$

\item
$X$ and $Y$ are absolutely continuous with regard to Lebesgue's measure.
Their densities are bounded and continuous on $\g R$.
\end{enumerate}
\end{Th}

\begin{Cor}
\label{densiteWDT}
For any $(\alpha,\beta )$, the distribution of $W_{(\alpha,\beta)}^{DT}$
admits a bounded and continuous density.
%is absolutely continuous with regard to Lebesgue's measure.
Its support is the whole real line $\g R$.
\end{Cor}

Corollary~\ref{densiteWDT} is an immediate consequence of Theorem~\ref{densiteXYDT} and of the decomposition property~\eqref{decompDT}.
The proof of Theorem~\ref{densiteXYDT}, that follows Liu's method, shows successively that
the distributions of $X$ and $Y$ are supported by the whole real line, that the characteristic functions $\varphi _X$ and $\varphi _Y$ reach the value $1$ only at the origin, that they tend to zero at $\pm\infty$ and finally that they are bounded above, in a neighbourhood of infinity, by a suitable power function so that an Fourier inversion theorem can apply, revealing the absolute continuity.

We first show a couple of lemmas, Lemma~\ref{supportWDT} being the first item of Theorem~\ref{densiteXYDT}.

\begin{Lem}
\label{supportWDT}
The support of both $X$ and $Y$ is the whole real line $\g R$.
\end{Lem}

{\sc Proof of Lemma~\ref{supportWDT}.\ }
We denote by $\Supp (X)$ and $\Supp (Y)$ the supports of $X$ and $Y$.
Since $\g EX>0$ and $\g EY<0$ (see Formula~\eqref{esperanceDT}), let $x>0$ and $y<0$ respectively belong to $\Supp (X)$ and $\Supp (Y)$.
Because of fixed point Equation~\eqref{pointFixeDiscret}, for any $v=(v_1,\dots ,v_{S+1})$ and
$w=(w_1,\dots ,w_{S+1})$ in $[0,1]$ such that $\sum _{1\leq k\leq S+1}v_k=\sum _{1\leq k\leq S+1}w_k=1$,
\begin{equation}
\label{eqPtFixeSupport}
\left( x\sum _{k=1}^{a+1}v_k^\sigma +y\sum _{k=a+2}^{S+1}v_k^\sigma ,
x\sum _{k=1}^{c}w_k^\sigma +y\sum _{k=c+1}^{S+1}w_k^\sigma\right)
\in\Supp (X)\times\Supp (Y) .
\end{equation}
We proceed in three steps.

Step 1:
{\it there exists $\varepsilon >0$ such that $[-\varepsilon ,\varepsilon ]\subseteq\Supp (X)\cap\Supp (Y)$}.

Indeed, apply Formula~\eqref{eqPtFixeSupport} for $v=w=(t,0,\dots ,0,1-t)$ where $t\in [0,1]$.
Then, $[y,x]=\{ t^\sigma x+(1-t)^\sigma y,~t\in [0,1]\}\subseteq\Supp (X)\cap\Supp (Y)$.
It suffices to take $\varepsilon =\min\{ x,-y\}$.

Step 2:
{\it there exists $\eta >0$ such that, for any $z\in\g R$},
$$
z\in\Supp (X)\cap\Supp (Y)\Longrightarrow \left( 1+\eta\right) z\in\Supp (X)\cap\Supp (Y).
$$
Indeed, apply Formula~\eqref{eqPtFixeSupport} for $v=\left( \frac 1{a+1},\dots ,\frac 1{a+1},0,\dots ,0\right)$ and $w=\left( 0,\dots ,0,\frac 1{d+1},\dots ,\frac 1{d+1}\right)$.
Then, $z\in\Supp (X)\cap\Supp (Y)$ implies that $\left( 1+a\right)^{1-\sigma}z\in\Supp (X)$ and $\left( 1+d\right)^{1-\sigma}z\in\Supp (Y)$.
It suffices to take $\eta =\min\{ \left( 1+a\right)^{1-\sigma}-1,\left( 1+d\right)^{1-\sigma}-1\}$.

Step 3:
the images of $[-\varepsilon ,\varepsilon ]$ by the iterates of the homothetic transformation $z\mapsto (1+\eta )z$ fill the whole real line.
\QED

\begin{Lem}
\label{phi<1}
For any $t\neq 0$, $\left| \varphi _X(t)\right| <1$ and $\left| \varphi _Y(t)\right| <1$.
\end{Lem}

{\sc Proof of Lemma~\ref{phi<1}.\ }
Of course, $\left| \varphi _X(t)\right|\leq1$ for any real number $t$.
Assume that $t_0\in\g R$ satisfies $\left| \varphi _X\left( t_0\right)\right| =1$.
Let $\theta _0\in\g R$ such that $\g E\left( e^{it_0X}\right) =e^{i\theta _0}$.
Then, almost surely, $e^{it_0X}=e^{i\theta _0}$ which is possible only if $t_0=0$ (and $\theta _0\in 2\pi\g Z$) since $\Supp (X)=\g R$.
Same proof for $Y$.
\QED

\begin{Lem}
\label{phi->0}
$\displaystyle\lim _{t\to\pm\infty}\varphi _X(t)=0$
and $\displaystyle\lim _{t\to\pm\infty}\varphi _Y(t)=0$
\end{Lem}

{\sc Proof of Lemma~\ref{phi->0}.\ }
This proof and the remainder of the argumentation on the absolute continuity rely on the following equalities that are consequences of the fixed point Equation~\eqref{pointFixeDiscret}:
for any $t\in\g R$,
\begin{equation}
\label{egalitePhis}
\left\{
\begin{array}{l}
\varphi _X(t)=\g E\Big(
\varphi _X\left( tV_1^\sigma\right)\dots\varphi _X\left( tV_{a+1}^\sigma\right)
\varphi _Y\left( tV_{a+2}^\sigma\right)\dots\varphi _Y\left( tV_{S+1}^\sigma\right)
\Big)
\\ \\
\varphi _Y(t)=\g E\Big(
\varphi _X\left( tV_1^\sigma\right)\dots\varphi _X\left( tV_{c}^\sigma\right)
\varphi _Y\left( tV_{c+1}^\sigma\right)\dots\varphi _Y\left( tV_{S+1}^\sigma\right)
\Big)
\end{array}
\right.
\end{equation}
where $V=\left( V_1,\dots ,V_{S+1}\right)$ is a Dirichlet distributed random vector with parameters $\left(\frac 1S,\dots ,\frac 1S\right)$.
These relations are obtained from fixed point Equation~\eqref{pointFixeDiscret} by conditioning with respect to~$V$.
In particular, since all $V_k$ are no zero with positive probability, Fatou's Lemma together with Equations~\eqref{egalitePhis} imply that
$$
\left\{
\begin{array}{l}
\limsup _{t\to\pm\infty}\left|\varphi _X(t)\right|
\leq\Big(\limsup _{t\to\pm\infty}\left|\varphi _X(t)\right|\Big) ^{a+1}
\\ \\
\limsup _{t\to\pm\infty}\left|\varphi _Y(t)\right|
\leq\Big(\limsup _{t\to\pm\infty}\left|\varphi _Y(t)\right|\Big) ^{d+1}.
\end{array}
\right.
$$
Consequently, since $a\geq 1$ and $d\geq 1$ (the urn is assumed to be large), $\limsup _{t\to\pm\infty}\left|\varphi _X(t)\right|\in\{ 0,1\}$ and the same holds for $\varphi _Y$.
It remains to show that $\limsup _{t\to\pm\infty}\left|\varphi _X(t)\right| =1$ is impossible to get the result for $X$
(a same argument for $Y$ applies as well).
The first Equation~\eqref{egalitePhis} implies that
$$
\left|\varphi _X(t)\right|\leq\g E\left|\varphi _X\left( V_1^\sigma t\right)\right| .
$$
The end of the proof relies on the following idea:
denoting by $U$ the uniform distribution on $[0,1]$, since $V_1^\sigma\egalLoi U^m\leq 1$ almost surely and $EV_k^\sigma =\frac 1{m+1}<1$ for any $k$, iterating this last inequality leads to $\limsup _{t\to\pm\infty}\left|\varphi _X(t)\right| =0$ which implies the final result.
The details, that are rather technical, can be almost literally adapted from Liu's proof of a result of the same kind.
See~\cite{Liu01}, Lemma 3.1, page 93.
On can also refer to~\cite{ChaLiuPouContinu} for a similar argument in a slightly different context.
\QED

\begin{Lem}
\label{majorationPhi}
For any $\rho '\in ]0,\frac 1m[$, when $t$ tends to $\pm\infty$, $\varphi _X(t)\in O\left( t^{-\rho '}\right)$ and $\varphi _Y(t)\in O\left( t^{-\rho '}\right)$.
\end{Lem}

{\sc Proof of Lemma~\ref{majorationPhi}.\ }
Let $\varepsilon >0$.
Let $T>0$ such that $\left|\varphi _X(t)\right|\leq\varepsilon$ and $\left|\varphi _Y(t)\right|\leq\varepsilon$ as soon as $|t|\geq T$;
the existence of $T$ is guaranteed by Lemma~\ref{phi->0}.
Then, because of~\eqref{egalitePhis}, for any $t\in\g R$, 
$$
\left|\varphi _X(t)\right|
\leq \varepsilon ^S\g E\left|\varphi _X\left( V_{S+1}^\sigma t\right)\right|
+\sum _{k=1}^S\g P\left( V_k^\sigma |t|\leq T\right).
$$
Since any $V_k^\sigma \egalLoi U^m$  where $U$ denotes the uniform distribution on $[0,1]$, this leads to
$$
\left|\varphi _X(t)\right|
\leq \varepsilon ^S\g E\left|\varphi _X\left( U^m t\right)\right|
+S\left(\frac T{|t|}\right)^{\frac 1m}
$$
for any $t\in\g R\setminus\{ 0\}$.
Now, for any $\rho\in ]0,1/m[$, $\g E\left( U^{-m\rho}\right) <\infty$ and the former inequality implies that there exists a positive constant $C$ such that for any nonzero $t$,
$$
\left|\varphi _X(t)\right|
\leq \varepsilon ^S\g E\left|\varphi _X\left( U^{m} t\right)\right|
+C\left(\frac 1{|t|}\right)^{\rho}.
$$
Thus, the random variables $X$ and $U$ satisfy the assumptions of the Gronwall type Lemma shown in~\cite{Liu01}, Lemme 3.2 page 93.
Using iterations of the former inequality, one gets for any $n\geq 1$,
$$
\left|\varphi _X(t)\right|
\leq\varepsilon ^{ns}\g E\left|\varphi _X\left( U_1^m\dots U_n^{m} t\right)\right|
+C|t|^{-\rho}\sum _{k=0}^{n-1}\Big(\varepsilon^S\g E\left( U^{-m\rho}\right)\Big) ^{k},
$$
which entails that
$\left|\varphi _X(t)\right|\leq C|t|^{-\rho}/\left( 1-\varepsilon^S\g E\left( U^{-m\rho}\right)\right)$
as soon as $\varepsilon$ is chosen in order that $1-\varepsilon^S\g E\left( U^{-m\rho}\right)>0$.
This implies the result.
A same argument is used for $\varphi _Y$.
\QED

{\sc End of the proof of Theorem~\ref{densiteXYDT}.\ }
Let $\rho\in ]0,\frac{a+1}m[$ and let $\rho '=\frac{\rho}{a+1}$.
Let $\kappa >0$ such that $\varphi _X(t)\leq \kappa |t|^{-\rho '}$ for any $t\neq 0$;
the existence of $\kappa$ is due to Lemma~\ref{majorationPhi}.
Applying~\eqref{egalitePhis}, one gets the successive inequalities
$$
\left|\varphi _X(t)\right|
\leq
\g E\left( \prod _{k=1}^{a+1}\left|\varphi _X\left( V_k^\sigma t\right)\right|\right)
\leq
\frac{\kappa ^{a+1}}{\left| t\right| ^\rho}\g E\left(\prod _{k=1}^{a+1}V_k^{-\sigma\rho '}\right)
$$
as soon as the last expectation is defined.
Since the random vector $V=\left( V_1,\dots ,V_{S+1}\right)$ is Dirichlet distributed with parameters $\left(\frac 1S,\dots ,\frac 1S\right)$, this expectation can be computed from the Stieltjes transform of $V$
(see the Appendix for the general form of joint moments, that can be extended to nonreal powers):
$$
\g E\left(\prod _{k=1}^{a+1}V_k^{-\sigma\rho '}\right)
=\frac{\Gamma\left( 1+\frac 1S\right)}{\Gamma\left( 1+\frac 1S-(a+1)\sigma\rho '\right)}
\left(\frac{\Gamma\left(\frac 1S-\sigma\rho '\right)}{\Gamma\left(\frac 1S\right)}\right) ^{a+1}
$$
is finite since $\sigma\rho '<\frac\sigma m=\frac 1S$ and $1+\frac 1S-(a+1)\sigma\rho '>1-\frac aS>0$.
Note that $a<S=a+b$ because the urn assumed to be is non triangular.
This proves \emph{(ii)}.
The same argument is used for the similar result on $\varphi _Y$ {\it (iii)}.

Since $\frac {a+1}m=\frac{a+1}{a-c}>1$, item {\it (ii)} implies that the Fourier transform $\varphi _X$ of the probability measure of $X$ is integrable.
This implies that $X$ admits a bounded continuous function as density.
The same result holds for $Y$.
\QED

\begin{Rem}
For the continuous time urn process, as shown in [ChaPouSah], the limit random variables $W_{(1,0)}^{CT}$ and  $W_{(0,1)}^{CT}$ admit densities as well.
These functions have been shown to be infinitely differentiable outside $0$, monotonic on $\g R_{<0}$ and $\g R_{>0}$, but not bounded around the origin.
The different behaviours of $W^{CT}$ and $W^{DT}$ have to be related to the martingale connexion~\eqref{connection}:
when the process starts with one ball, the density of $\xi ^{\sigma}$ is not bounded at $0$ since $\xi$ is Gamma distributed, with parameter $\frac 1S$.
\end{Rem}

%%%%%%%%%%%%%%%%%%%%%%%%%%%%%%%%%
\section{Smoothing transforms}
\label{sec:contractions}
%%%%%%%%%%%%%%%%%%%%%

This section is devoted to the existence and the unicity of solutions of the distributional systems~\eqref{pointFixeDiscret} and~\eqref{pointFixeContinu}.
Notice that existence and unicity of solutions of the discrete-time system~\eqref{pointFixeDiscret} could be deduced from the general result in Neininger-R\"uschendorf \cite{NeiRusaap}, nevertheless we give hereunder a rapid and autonomous proof of Theorem \ref{th:solutionsSystemeDiscret}, in order to make explicit the contraction method in the case of large P\'olya urns. The proof is reminiscent of the one in Fill-Kapur \cite{FillKapur}.

When $A$ is a real number, let $\rond M_2\left( A\right)$ be the space of probability distributions on $\g R$ that have $A$ as expectation and a finite second moment, endowed with a complete metric space structure by the Wasserstein distance.
Note first that when $X$ and $Y$ are solutions of ~\eqref{pointFixeDiscret} or~\eqref{pointFixeContinu} that have respectively $B$ and $C$ as expectations, then $cB+bC=0$ (elementary computation).
In Theorems~\ref{th:solutionsSystemeDiscret} and~\ref{th:solutionsSystemeContinu}, we prove that when $B$ and $C$ are two real numbers that satisfy $cB+bC=0$, the systems~\eqref{pointFixeDiscret} and~\eqref{pointFixeContinu} both have a unique solution in the product metric space $\rond M_2\left( B\right)\times\rond M_2\left( C\right)$.
To do so, we use the Banach contraction method.

Since $\left(\g EX,\g EY\right)$ is proportional to $\left( b,-c\right)$ in both continuous time and discrete time urn processes (Formulae~\eqref{esperanceDT} and~\eqref{esperanceCT}), this result shows that the systems~\eqref{pointFixeDiscret} and~\eqref{pointFixeContinu} characterize the limit distributions $W_{(1,0)}^{DT}$ and $W_{(0,1)}^{DT}$ on one hand, $W_{(1,0)}^{CT}$ and $W_{(0,1)}^{CT}$ on the other hand.

%%%%%%%%%%%%%%%%%%%%%%%%%%
\subsection{The Wasserstein distance}

Let $A\in\g R$.
The Wasserstein distance on $\rond M_2\left( A\right)$ is defined as follows:
$$
d_W\left(\mu _1,\mu _2\right) = \min _{\left( X_1, X_2\right)} \Big( \g E\left( X_1-X_2\right) ^2\Big) ^{1/2}
$$
where the minimum is taken over random vectors $\left( X_1,X_2\right)$ on $\g R^2$ having respective marginal distributions $\mu _1$ and $\mu _2$ ;
the minimum is attained  by the Kantorovich-Rubinstein Theorem.
With this distance, $\rond M_2\left( A\right)$ is a complete metric space (see for instance Dudley \cite{Dudley}).

Let $(B,C)\in\g R^2$.
The product space $\rond M_2\left( B\right)\times\rond M_2\left( C\right)$ is equipped with the product metric, defined (for example) by the distance
$$
d\Big( \left( \mu _1,\nu _1\right), \left( \mu _2,\nu _2\right)\Big)
= \max\Big\{ d_W\left(\mu_1,\mu _2\right) ,d_W\left(\nu_1,\nu _2\right)\Big\} .
$$
Of course, this product remains a complete metric space.

%%%%%%%%%%%%%%%%%%%
\subsection{Contraction method in discrete time}
\label{sec:contractionDiscrete}

Let us recall the fixed point system~\eqref{pointFixeDiscret} satisfied by $(X^{DT},Y^{DT})$, the elementary limits of a large two-color discrete time P\'olya urn process:
\begin{equation*}
\left\{
\begin{array}{l}
\displaystyle
X\egalLoi\sum_{k=1}^{a+1} V_k^{\sigma} X^{(k)} + \sum_{k=a+2}^{S+1} V_k^{\sigma} Y^{(k)}\\ \\
\displaystyle
Y\egalLoi\sum_{k=1}^{c}V_k^{\sigma} X^{(k)} + \sum_{k=c+1}^{S+1}V_k^{\sigma} Y^{(k)}.
\end{array}
\right.
\end{equation*}

Let $\rond M_2$ be the space of square-integrable probability measures on $\g R$.
When $(B,C)\in\g R^2$, let $K_1$ be the function defined on $\rond M_2\left( B\right)\times\rond M_2\left( C\right)$ by:
$$
\begin{array}{rccl}
K_1:
&\rond M_2\left( B\right)\times\rond M_2\left( C\right)&\longrightarrow&\rond M_2\\
&(\mu ,\nu)&\longmapsto
&\displaystyle\rond L\left(\sum_{k=1}^{a+1} V_k^{\sigma} X^{(k)} + \sum_{k=a+2}^{S+1} V_k^{\sigma} Y^{(k)}\right)
\end{array}
$$
where $X^{(1)},\dots ,X^{(a+1)}$ are $\mu$-distributed random variables, 
$Y^{(a+2)},\dots ,Y^{(S+1)}$ are $\nu$-distributed random variables,
$V=\left( V_1,\dots ,V_{S+1}\right)$ is a Dirichlet-distributed random vector with parameter $\left(\frac 1S,\dots ,\frac 1S\right)$, the $X^{(k)}$, $Y^{(k)}$ and $V$ being all independent of each other.
Similarly, let $K_2$ be defined by
$$
\begin{array}{rccl}
K_2:
&\rond M_2\left( B\right)\times\rond M_2\left( C\right)&\longrightarrow&\rond M_2\\
&(\mu ,\nu)&\longmapsto
&\displaystyle\rond L\left(\sum_{k=1}^{c} V_k^{\sigma} X^{(k)} + \sum_{k=c+1}^{S+1} V_k^{\sigma} Y^{(k)}\right) .
\end{array}
$$
A simple computation shows that if $(\mu ,\nu )\in\rond M_2\left( B\right)\times\rond M_2\left( C\right)$, then
$$\g EK_1(\mu ,\nu )=\frac{(a+1)B+bC}{m+1}$$
and
$$\g EK_2(\mu ,\nu )=\frac{cB+(d+1)C}{m+1},$$
so that, since $m=a-c=d-b$, the relation $cB+bC=0$ is a sufficient and necessary condition for the product function $\left( K_1,K_2\right)$ to range $\rond M_2\left( B\right)\times\rond M_2\left( C\right)$ into itself.

\begin{Lem}
\label{th:contractionDT}
Let $B$ and $C$ be real numbers that satisfy $cB+bC=0$.
Then, the smoothing transform
$$
\begin{array}{rccl}
K:&\rond M_2\left( B\right)\times\rond M_2\left( C\right)
&\longrightarrow&\rond M_2\left( B\right)\times\rond M_2\left( C\right)\\
&(\mu ,\nu )&\longmapsto&\Big( K_1(\mu ,\nu ),K_2(\mu ,\nu )\Big)
\end{array}
$$
is $\sqrt{\frac{S+1}{2m+1}}$-Lipschitz.
In particular, it is a contraction.
\end{Lem}

\begin{Th}
\label{th:solutionsSystemeDiscret}
\begin{enumerate}[(i)]
\item
When $B$ and $C$ are real numbers that satisfy $cB+bC=0$, System~\eqref{pointFixeDiscret} has a unique solution in $\rond M_2\left( B\right)\times\rond M_2\left( C\right)$.
\item
The pair $\left( X^{DT},Y^{DT}\right)$ is the unique solution of the distributional System~\eqref{pointFixeDiscret} having
$\left( \frac{\Gamma\left(\frac 1S\right)}{\Gamma\left(\frac {m+1}S\right)}\frac{b}{S},-\frac{\Gamma\left(\frac 1S\right)}{\Gamma\left(\frac {m+1}S\right)}\frac{c}{S}\right)$
as expectation and a finite second moment.
\end{enumerate}
\end{Th}

Theorem~\ref{th:solutionsSystemeDiscret} is a direct consequence of Lemma~\ref{th:contractionDT} and of Banach's fixed point theorem.

\vskip 5pt
{\sc Proof of Lemma~\ref{th:contractionDT}.}
Let $(\mu _1,\nu _1)$ and $(\mu _2,\nu _2)$ in
$\rond M_2\left( B\right)\times\rond M_2\left( C\right)$.
Let $V=\left( V_1,\dots ,V_{S+1}\right)$ be a Dirichlet random vector with parameter $\left(\frac 1S,\dots ,\frac 1S\right)$.
Let $X_1^{(1)},\dots ,X_1^{(a+1)}$ be $\mu _1$-distributed random variables, $Y_1^{(a+2)},\dots ,Y_1^{(S+1)}$ be $\nu _1$-distributed random variables,
$X_2^{(1)},\dots ,X_1^{(c)}$ be $\mu _2$-distributed random variables and
$Y_2^{(c+1)},\dots ,Y_2^{(S+1)}$ be $\nu _2$-distributed random variables, all of them being independent and independent of $V$.
Then,
$$
\begin{array}{l}
\displaystyle
d_W\Big( K_1\left(\mu _1,\nu _1\right) ,K_1\left(\mu _2,\nu _2\right)\Big) ^2
\leq
\left\|
\sum_{k=1}^{a+1} V_k^{\sigma}\left( X_1^{(k)}-X_2^{(k)}\right)
+\sum_{k=a+2}^{S+1} V_k^{\sigma}\left( Y_1^{(k)}-Y_2^{(k)}\right)
\right\|_2^2
\\ \displaystyle
\hskip 100pt=\Var\left[
\sum_{k=1}^{a+1} V_k^{\sigma}\left( X_1^{(k)}-X_2^{(k)}\right)
+\sum_{k=a+2}^{S+1} V_k^{\sigma}\left( Y_1^{(k)}-Y_2^{(k)}\right)
\right]
\\ \displaystyle
\hskip 100pt=\g E\Var\left(\left.
\sum_{k=1}^{a+1} V_k^{\sigma}\left( X_1^{(k)}-X_2^{(k)}\right)
+\sum_{k=a+2}^{S+1} V_k^{\sigma}\left( Y_1^{(k)}-Y_2^{(k)}\right)
\right| V\right)
\\ \displaystyle
\hskip 110pt+\Var\g E\left(\left.
\sum_{k=1}^{a+1} V_k^{\sigma}\left( X_1^{(k)}-X_2^{(k)}\right)
+\sum_{k=a+2}^{S+1} V_k^{\sigma}\left( Y_1^{(k)}-Y_2^{(k)}\right)
\right| V\right)
\end{array}
$$
thanks to the law of total variance.
Since $V=\left( V_1,\dots ,V_{S+1}\right)$ is independent of the $X_j^{(k)}$ and of the $Y_j^{(k)}$, one gets
\begin{align*}
d_W\Big( K_1\left(\mu _1,\nu _1\right) ,K_1\left(\mu _2,\nu _2\right)\Big) ^2
& \leq
\sum_{k=1}^{a+1}\g EV_k^{2\sigma}\Var\left( X_1^{(k)}-X_2^{(k)}\right)
+\sum_{k=a+2}^{S+1}\g EV_k^{2\sigma}\Var\left( Y_1^{(k)}-Y_2^{(k)}\right)
\\
&\leq
\Var \left( X_1^{(1)}-X_2^{(1)}\right)
\sum_{k=1}^{a+1}\g EV_k^{2\sigma}
+\Var \left( Y_1^{(1)}-Y_2^{(1)}\right)
\sum_{k=a+2}^{S+1}\g EV_k^{2\sigma}
\\
&=
\frac{a+1}{2m+1}\left\| X_1^{(1)}-X_2^{(1)}\right\| _2^2 
+\frac{b}{2m+1}\left\| Y_1^{(1)}-Y_2^{(1)}\right\| _2^2.
\end{align*}
Since the inequality holds for any random variables $X_1^{(1)}$, $X_2^{(1)}$, $Y_1^{(1)}$ and $Y_2^{(1)}$ having respective distributions $\mu _1$, $\mu _2$, $\nu _1$ and $\nu _2$, this leads to
\begin{align*}
d_W\Big( K_1\left(\mu _1,\nu _1\right) ,K_1\left(\mu _2,\nu _2\right)\Big) ^2
&\leq
\frac{a+1}{2m+1}d_W\left(\mu _1,\mu _2\right) ^2
+\frac{b}{2m+1}d_W\left(\nu _1,\nu _2\right) ^2\\
&\leq
\frac{S+1}{2m+1}d\Big( \left( \mu _1,\nu _1\right), \left( \mu _2,\nu _2\right)\Big) ^2.
\end{align*}
A very similar computation shows that
$$
d_W\Big( K_2\left(\mu _1,\nu _1\right) ,K_2\left(\mu _2,\nu _2\right)\Big) ^2
\leq
\frac{S+1}{2m+1}d\Big( \left( \mu _1,\nu _1\right), \left( \mu _2,\nu _2\right)\Big) ^2,
$$
so that, finally,
$$
d\Big( K\left(\mu _1,\nu _1\right) ,K\left(\mu _2,\nu _2\right)\Big) ^2
\leq
\frac{S+1}{2m+1}d\Big( \left( \mu _1,\nu _1\right), \left( \mu _2,\nu _2\right)\Big) ^2
$$
making the proof complete.
Note that the assumption $\sigma =\frac mS>\frac 12$ guarantees that the Lipschitz constant is in $]0,1[$.
\QED

%%%%%%%%%%%%%%%%%%%%%%%%%%%%%%%%%
\subsection{Contraction method in continuous time}

In continuous time, the laws of $X^{CT}$ and $Y^{CT}$ are solutions of the following system (cf.~\eqref{pointFixeContinu}):
\begin{equation*}
\left\{
\begin{array}{l}
\displaystyle
X\egalLoi U^{m}\bigg( \sum_{k=1}^{a+1}X^{(k)}+\sum_{k=a+2}^{S+1}Y^{(k)}\bigg)\\ \\
\displaystyle Y\egalLoi U^{m}\bigg(\sum_{k=1}^{c}X^{(k)}+\sum_{k=c+1}^{S+1}Y^{(k)} \bigg),
\end{array}
\right.
\end{equation*}
The following theorem, which is the continuous time version of Theorem~\ref{th:solutionsSystemeDiscret}, can be proved by adapting the arguments of Theorem~\ref{th:solutionsSystemeDiscret}.
Details are left to the reader.

\begin{Th}
\label{th:solutionsSystemeContinu}
\begin{enumerate}[(i)]
\item
When $B$ and $C$ are real numbers that satisfy $cB+bC=0$, System~\eqref{pointFixeContinu} has a unique solution in $\rond M_2\left( B\right)\times\rond M_2\left( C\right)$.
\item
The pair $\left( X^{CT},Y^{CT}\right)$ is the unique solution of the distributional System~\eqref{pointFixeContinu} having
$\left(\frac{b}{S},-\frac{c}{S}\right)$
as expectation and a finite second moment.
\end{enumerate}
\end{Th}

%%%%%%%%%%%%%%%%%%%%%%%%%%%%%%%%%%
\section{Moments}
\label{sec:moments}
%%%%%%%%%%%%%%%%%%%%%%%%%%%%%

This section is devoted to the asymptotics of the moments of the limit variables $W^{DT}$ and $W^{CT}$. We shall see that they are big but not too much. Observe first that the connexion~\eqref{connection} allows us to study only one of the two cases among discrete or continuous case. We chose to focus on the continuous case, since the fixed point equation system is slightly easier to deal with. Let us recall here system \eqref{pointFixeContinu}.

\begin{equation*}
%\label{eqpointfixeBis}
\left\{
\begin{array}{l}
\displaystyle
X\egalLoi U^{m}\bigg( \sum_{k=1}^{a+1}X^{(k)}+\sum_{k=a+2}^{S+1}Y^{(k)}\bigg)\\ \\
\displaystyle Y\egalLoi U^{m}\bigg(\sum_{k=1}^{c}X^{(k)}+\sum_{k=c+1}^{S+1}Y^{(k)} \bigg),
\end{array}
\right.
\end{equation*}
where $U$ is uniform on $[0,1]$,
where $X$,  $X^{(k)}$ and $Y$,  $Y^{(k)}$ are respective copies of $X^{CT}$ and $Y^{CT}$, all being independent of
each other and of $U$.

\medskip

Up to now, what is known about the size of these moments is contained in \cite{ChaPouSah} where it is proved that the radius of convergence of the Laplace series of a non trivial square integrable solution of~\eqref{pointFixeContinu} is equal to zero.
Consequently, by the Hadamard formula for the radius of convergence, 
$$
\limsup_p \left( \frac{\g E |X|^p}{p!}\right)^{\frac 1p} = +\infty .
$$
In otherwords, for any constant $C$, for $p$ large enough,
$$
C^p \leq \frac{\g E |X|^p}{p!} .
$$
The following lemma gives an upperbound for $\frac{\g E |X|^p}{p!}$.  It is the argument leading to Theorem~\ref{detMoment} where it is proved that the law of $X$ is determined by its moments. 

\begin{Lem}
\label{borne}
If $X$ and $Y$ are integrable solutions of ~\eqref{pointFixeContinu}, they admit absolute moments of all orders $p\geq 1$ and the sequences $\displaystyle\left(\frac{\g E |X|^p}{p!\log ^pp}\right)^{\frac 1p}$ and $\displaystyle\left(\frac{\g E |Y|^p}{p!\log ^pp}\right)^{\frac 1p}$ are bounded.
\end{Lem}

\pff
%, taking advantage of system \eqref{eqpointfixeBis}. 
Let $\varphi(p):= \log ^p (p+2)$ and define
$$
u_p:=\displaystyle\frac{\E |X|^p}{p!\varphi(p)}
{\rm ~~and~~}
v_p:= \displaystyle\frac{\E |Y|^p}{p!\varphi(p)}.
$$
We show by induction on $p\geq 1$ that $\left(\frac{\g E |X|^p}{p!\varphi (p)}\right)^{\frac 1p}$ and $\left(\frac{\g E |Y|^p}{p!\varphi (p)}\right)^{\frac 1p}$ are finite and define bounded sequences. Notice that a similar technique is used in Kahane-Peyri\`ere \cite{KahPey}.
Take the power $p$ in the first equation notice that $\E U^{mp} = \displaystyle\frac 1{mp+1}$, and isolate the two extreme terms. One gets (remember $S+1 = a+1+b$)
\begin{equation*}
\begin{array}{ll}
\E |X|^p &\displaystyle\leq \frac 1{mp+1} \Bigg(  (a+1)\E |X|^p + b\E |Y|^p\\
&
\hskip 20pt +\displaystyle\summ{p_1+ \dots + p_{S+1} = p}{ p_j\leq p-1} \ \frac{p!}{p_1!\dots p_{S+1}!}\E |X|^{p_1} \dots \E |X|^{p_{a+1}}\E |Y|^{p_{a+2}}\dots \E |Y|^{p_{S+1}}\Bigg)
\end{array}
\end{equation*}
or also
$$
(mp-a)\E |X|^p \leq  b \E |Y|^p +  \summ{p_1+ \dots + p_{S+1} = p}{ p_j\leq p-1} \   \frac{p!}{p_1!\dots p_{S+1}!}\E |X|^{p_1} \dots \E |X|^{p_{a+1}}\E |Y|^{p_{a+2}}\dots \E |Y|^{p_{S+1}}  .
$$
An analog inequality holds for $\E |Y|^p$, leading to the system
\begin{equation}
\label{inegaliteMoments}
\left\{
\begin{array}{l}
\displaystyle
(mp-a) u_p \leq bv_p + \summ{p_1+ \dots + p_{S+1} = p}{ p_j\leq p-1} \  u_{p_1} \dots u_{p_{a+1}}v_{p_{a+2}}\dots v_{p_{S+1}} \frac{\varphi(p_1)\dots \varphi(p_{S+1})}{\varphi(p)}\\ \\
\displaystyle (mp-d) v_p \leq cu_p + \summ{p_1+ \dots + p_{S+1} = p}{ p_j\leq p-1} \   u_{p_1} \dots u_{p_{c}}v_{p_{c+1}}\dots v_{p_{S+1}} \frac{\varphi(p_1)\dots \varphi(p_{S+1})}{\varphi(p)}.
\end{array}
\right.
\end{equation}
Since the eigenvalues of the matrix $R=\begin{pmatrix}a&b\\c&d\end{pmatrix}$ are $m$ and $S$ and since $2m>S$ (the urn is assumed to be large), all matrices $mpI_2-R$ $(p\geq 2)$ are invertible so that System~\eqref{inegaliteMoments} implies by induction on $p$ that solutions $X$ and $Y$ of System~\eqref{pointFixeContinu} admit absolute moments of all orders as soon as they are integrable.

Let $p_0$ be the smallest positive integer such that for any $p\geq p_0$,
$$
\frac{m(p-1)}{(mp-a)(mp-d)-bc} \Big(1+8 \log\left( p+2\right)\Big)^{S+1} \leq 1.
$$
Such a $p_0$ exists since the left handside goes to $0$ when $p$ goes to $+\infty$. Denote
$$
A:= \max _{1\leq q\leq p_0}
\left\{\left( u_q\right) ^{\frac 1q}, \left( v_q\right) ^{\frac 1q}\right\}.
$$
Assume by induction on $p\geq p_0 +1$ that for every $q\leq p-1, (u_q)^{\frac 1q} \leq A$ and $(v_q)^{\frac 1q} \leq A$.
Then, 
\begin{equation*}
%\label{dislocationProjetee}
\left\{
\begin{array}{l}
\displaystyle
(mp-a) u_p \leq bv_p + A^p\summ{p_1+ \dots + p_{S+1} = p}{ p_j\leq p-1} \   \frac{\varphi(p_1)\dots \varphi(p_{S+1})}{\varphi(p)}\\ \\
\displaystyle (mp-d) v_p \leq cu_p + A^p\summ{p_1+ \dots + p_{S+1} = p}{ p_j\leq p-1} \   \frac{\varphi(p_1)\dots \varphi(p_{S+1})}{\varphi(p)}.
\end{array}
\right.
\end{equation*}
Let
\be
\label{Phi}
\Phi(p):= \summ{p_1+ \dots + p_{S+1} = p}{ p_j\leq p-1} \   \frac{\varphi(p_1)\dots \varphi(p_{S+1})}{\varphi(p)} 
\ee
so that
\begin{equation*}
\left\{
\begin{array}{l}
\displaystyle
(mp-a) u_p \leq bv_p + A^p\Phi(p)\\ \\
(mp-d) v_p \leq cu_p + A^p\Phi(p)
\end{array}
\right.
\end{equation*}
which implies
$$
u_p \leq \frac{m(p-1)}{(mp-a)(mp-d)-bc}A^p \Phi(p)
$$
and the same inequality for $v_p$ as well.
Admit for a while the following lemma.

\begin{Lem}
\label{lemmePhi}
%For $\varphi(p):= \log ^p (p+2)$ and $\Phi$ defined by (\ref{Phi}), 
For every $p\geq 2$,
$\Phi (p)\leq \Big(1+8  \log\left( p+2\right)\Big)^{S+1}.$
\end{Lem}
Consequently 
$$
u_p\leq \frac{m(p-1)}{(mp-a)(mp-d)-bc}A^p\Big( 1+8 \log\left( p+2\right)\Big) ^{S+1}.
$$
%Since for all $p\geq p_0$, $\displaystyle \frac{m(p-1)}{(mp-a)(mp-d)-bc}\left(1+8  \log (p+2)\right)^{S+1}\leq 1$,
By definition of $p_0$, this implies that $(u_p)^{\frac 1p}\leq A$ and the recurrence holds. \QED

\vskip 5pt
{\sc Proof of Lemma \ref{lemmePhi}}.
The definitions of $\varphi$ and $\Phi$ imply directly that
\begin{eqnarray*}
\Phi(p)&=
&\summ{p_1+ \dots + p_{S+1} = p}{ p_j\leq p-1}
\frac{\log ^{p_1}\left( p_1+2\right)
\dots
\log ^{p_{S+1}}\left( p_{S+1}+2\right)}{\log ^{p}\left( p+2\right)}\\
&=
&\summ{p_1+ \dots + p_{S+1} = p}{ p_j\leq p-1}
\left( 1+\frac{\log \left( 1-\frac{p-p_1}{p+2}\right)}{\log\left( p+2\right)}\right)^{p_1}
\cdots
\left( 1+\frac{\log \left( 1-\frac{p-p_{S+1}}{p+2}\right)}{\log\left( p+2\right)}\right)^{p_{S+1}}.
\end{eqnarray*}
Using $\log(1-u)\leq -u$ for all $u<1$ leads to
$$
\Phi (p)\leq
\summ{p_1+ \dots + p_{S+1} = p}{ p_j\leq p-1}
\left( 1-\frac{p-p_1}{(p+2)\log\left( p+2\right)}\right)^{p_1}
\cdots
\left( 1-\frac{p-p_{S+1}}{(p+2)\log\left( p+2\right)}\right)^{p_{S+1}}
$$
which can be written with an exponential to get, using again $\log(1-u)\leq -u$:
$$
\Phi (p)\leq
\summ{p_1+ \dots + p_{S+1} = p}{ p_j\leq p-1}
\exp\left\{
-\displaystyle\frac{p^2}{(p+2)\log\left( p+2\right)} \displaystyle\sum_{j=1}^{S+1} \frac{p_j}{p} \left( 1- \frac{p_j}{p}\right)
\right\}
$$
Let 
$\psi_p(x):=\exp\left( -\displaystyle\frac{p^2}{(p+2)\log\left( p+2\right)} x(1-x)\right)$, so that
\begin{eqnarray*}
\Phi(p) & \leq & \summ{p_1+ \dots + p_{S+1} = p}{ p_j\leq p-1}  \psi_p\left(\frac{p_1}{p} \right)\dots \psi_p\left(\frac{p_{S+1}}{p}\right )\\
& \leq &  \sum_{0\leq p_1, \dots , p_{S+1} \leq  p-1}\psi_p\left(\frac{p_1}{p} \right)\dots \psi_p\left(\frac{p_{S+1}}{p}\right )\\
&=& \left( \sum_{k=0}^{p-1} \psi_p\left(\frac{k}{p} \right)\right)^{S+1}.
\end{eqnarray*}
Elementary calculations lead then to
$$
\sum_{k=0}^{p-1} \psi_p\left(\frac{k}{p} \right) \leq 1+p \int_0^{1} \psi_p(t) dt
$$
and  for any $\alpha >0$
$$
\int_0^{1} \exp\left( -\alpha x(1-x)\right)  dt \leq \frac 4{\alpha}
$$
so that
$$
\int_0^{1} \psi_p(t) dt \leq 4 \frac{(p+2)\log\left( p+2\right)}{p^2}
$$
and the lemma holds.
\QED

The upperbound on the moments, obtained in Lemma \ref{borne} leads to the following theorem.

\begin{Th}
\label{detMoment}
Let $X$ and $Y$ be integrable solutions of any fixed point equation~\eqref{pointFixeDiscret} or~\eqref{pointFixeContinu}.
Then, $X$ and $Y$ admit absolute moments of all orders $p\geq 1$ and the probability distributions of $|X|$,  $|Y|$,  $X$ and $Y$  are determined by their moments.
% (though their Laplace series have a radius of convergence equal to $0$).
\end{Th}

\pff
By Lemma \ref{borne}, if $X$ and $Y$ are integrable solutions of~\eqref{pointFixeContinu}, they admit moments of all orders and,
when $p$ is large enough,
\begin{equation}
\label{inequality}
\left(\g E |X|^p \right)^{-\frac 1p} \geq C\frac{(p!)^{-\frac 1p}}{\log p}.
%\frac{\g E |X|^p}{(p!)} \leq C^p( \log p)^p
\end{equation}
Besides, by Stirling's formula, when $p$ tends to infinity,
$$
\frac{(p!)^{-\frac 1p}}{\log p} \sim \frac{e}{p\log p}
%\left(\g E |X|^p \right)^{-\frac 1p} \geq \frac{(p!)^{-\frac 1p}}{\log p} \sim \frac{e}{p\log p}
$$
which is the general term of a Bertrand divergent series. The Carleman's criterion applies, implying that $X$ and $Y$ are moment determined.

If $X$ and $Y$ are integrable solutions of~\eqref{pointFixeDiscret} and if $\xi$ is an independent $Gamma\left(\frac 1S\right)$-distributed random variable, then, thanks to Proposition~\ref{connexionSystems}, $\xi ^\sigma X$ and $\xi ^\sigma Y$ are integrable solutions of~\eqref{pointFixeContinu} so that they both satisfy Carleman's criterion.
This implies that $X$ and $Y$ are moment determined as well.
\QED

\begin{Cor}
For any initial composition $(\alpha ,\beta )$, the limit laws $W^{DT}_{(\alpha ,\beta )}$ and $W^{CT}_{(\alpha ,\beta )}$ of a large P\'olya urn process are determined by their moments.
%The limit laws $W^{DT}_{(1,0)}$, $W^{DT}_{(0,1)}$, $W^{CT}_{(1,0)}$, $W^{CT}_{(0,1)}$ of a large P\'olya urn process are determined by their moments. The limit laws $W^{DT}_{(\alpha ,\beta )}$ and $W^{CT}_{(\alpha ,\beta )}$ as well. 
\end{Cor}

\pff
For elementary initial compositions $(1,0)$ or $(0,1)$, the result is a direct consequence of Theorems~\ref{th:solutionsSystemeDiscret}, \ref{th:solutionsSystemeContinu} and~\ref{detMoment}.
For a general initial composition $(\alpha ,\beta )$ in continuous time, notice that decomposition Formula~\eqref{decompCT} implies that
$$
|| W^{CT}_{(\alpha ,\beta )}||_p \leq \alpha ||W^{CT}_{(1,0)}||_p + \beta ||W^{CT}_{(0,1)}||_p.
$$
Since $W^{CT}_{(1,0)}$ and $W^{CT}_{(0,1)}$ satisfy \eqref{inequality}, $W^{CT}_{(\alpha ,\beta )}$ satisfies Carleman's criterion;
it is thus determined by its moments.
The same arguments hold in discrete time, using decomposition Formula~\eqref{decompDT}.
\QED

%%%%%%%%%%%%%%%%%%
\section{Appendix: P\'olya urns and Dirichlet distribution}
\label{sec:appendix}
%%%%%%%%%%%%%%%%%%

In this section, we deal with results that belong to the ``folklore'':
they are not new neither very difficult, but are nowhere properly gathered, to the best of our knowledge.
Proposition~\ref{proDirichlet} goes back to Athreya~\cite{Athreya69} with different names and a different proof.
It is partially given in Blackwell and Kendall \cite{BlaKen} for $S=1$ and starting from one ball of each color. 
The moment method is evocated in Johnson and Kotz book~\cite{JK}.
We detail here a proof to make our paper self-contained. 

%%%%%%%%%%%%%%%%
\subsection{Dirichlet distributions}
\label{Dirichlet}

This section gathers some well known facts on Dirichlet distributions.
Besides, we fix notations we use in the sequel.

Let $d\geq 2$ be a natural integer.
Let $\Sigma$ be the $(d-1)$-dimensional simplex
$$
\Sigma =\left\{ (x_1,\dots x_{d})\in [0,1]^{d},~\sum _{k=1}^{d}x_k=1\right\} .
$$
The following formula is a generalization of the definition of Euler's B\^eta function:
let $\left(\nu _1,\dots ,\nu _d\right)$ be positive real numbers.
Then,
\begin{equation}
\label{betaGen}
\begin{array}{rcl}
\displaystyle
\int _\Sigma\left[\prod_{k=1}^{d}x_k^{\nu _k-1}\right]d\Sigma \left(x_1,\dots ,x_{d}\right)
&=&\displaystyle
\frac{\Gamma (\nu _1)\dots\Gamma (\nu _{d})}{\Gamma (\nu _1+\dots +\nu _{d})}
\end{array}
\end{equation}
where $d\Sigma$ denotes the positive measure on the simplex $\Sigma$, defined by
$$
\begin{array}{c}
f\left( x_1,\dots ,x_{d}\right) d\Sigma\left( x_1,\dots ,x_{d}\right) \hskip 200pt\\
\hskip 80pt
=f\left( x_1,\dots ,x_{d-1},1-\sum _{k=1}^{d-1}x_k\right)
\1_{\left\{x\in\left[ 0,1\right] ^{d-1},~\sum _{k=1}^{d-1}x_k\leq 1\right\}}
dx_1\dots dx_{d-1}
\end{array}
$$
for any continuous function $f$ defined on $\Sigma$.

\vskip 10pt
By means of this formula, one defines usually the \emph{Dirichlet distribution} with parameters
$\left(\nu _1,\dots ,\nu _d\right)$, denoted by $Dirichlet\left(\nu _1,\dots ,\nu _d\right)$,
whose density on $\Sigma$ is given by
$$
\frac{\Gamma (\nu _1+\dots +\nu _{d})}{\Gamma (\nu _1)\dots\Gamma (\nu _{d})}
\left[\prod _{k=1}^{d}x_k^{\nu _k-1}\right]
d\Sigma \left(x_1,\dots ,x_{d}\right).
$$
In particular, if $D=(D_1,\dots ,D_d)$ is a $d$-dimensional random vector which is
Dirichlet-distributed with parameters $\left(\nu _1,\dots ,\nu _d\right)$, then, for any
$p=(p_1,\dots ,p_{d})\in\g N^{d}$, the (joint) moment of order $p$ of $D$ is
$$
\g E\left( D^p\right)
=\g E\left( D_1^{p_1}\dots D_d^{p_d}\right)
=\frac{\Gamma\left(\nu\right)}{\Gamma\left(\nu +|p|\right)}
\prod _{k=1}^{d}\frac{\Gamma\left( \nu _k+p_k\right)}{\Gamma \left(\nu _k\right)}
$$
where $\nu =\sum _{k=1}^d\nu _k$ and $|p|=\sum _{k=1}^dp_k$.

Finally, a computation of same kind shows that the $[0,1]$-valued random variable $D_k$, which is
the $k$-th marginal distribution of $D$, is $Beta\left( \nu _k, \nu-\nu _k\right)$-distributed
\emph{i.e.} admits the density
$$
\frac 1{B\left( \nu _k, \nu-\nu _k\right)}t^{\nu _k-1}\left( 1-t\right) ^{\nu -\nu _k-1}\1 _{[0,1]}dt.
$$

Note that computing asymptotics of such moments when $p$ tends to infinity by Stirling's formula
leads to show that a Dirichlet distribution is determined by its moments.
An alternative description of a Dirichlet distribution can be made by considering a sequence $\left( G_1,\dots ,G_d\right)$ of \emph{Gamma}-distributed random variables conditioned to the relation $\sum _{k=1}^dG_k=1$.

%%%%%%%%%%%%%%%%%%%%%%%%%%%%%%%%%
\subsection{Original/diagonal P\'olya urns}
%%%%%%%%%%%%%%%%%%%%%%%%%%%%%%%%%%

\begin{Prop}
\label{proDirichlet}
Let $d\geq 2$ and $S\geq 1$ be integers.
Let also $\left( \alpha _1,\dots ,\alpha _d\right)\in\g N^d\setminus \{ 0\}$.
%Denote $\alpha =\sum _{k=1}^d\alpha _k\geq 1$.
Let $\left( P_n\right) _{n\geq 0}$ be the $d$-color P\'olya urn random process having $SI_d$ as replacement matrix and $\left( \alpha _1,\dots ,\alpha _d\right)$ as initial composition.
Then, almost surely and in any ${\rm L}^t$, $t\geq 1$,
$$
\frac {P_n}{nS}\ \conv V
%\frac {P_n}{\alpha +nS}\ \conv V
$$
where $V$ is a $d$-dimensional Dirichlet-distributed random vector, with parameters
$(\frac{\alpha _1}{S},\dots ,\frac{\alpha _d}{S})$.
\end{Prop}

\begin{Rem}
For any $k\in\{ 1,\dots  ,d\}$, the $k$-th coordinate of $V$ is
%$Beta (\frac{\alpha_k}{S},\frac{\sum_{j=1}^d \alpha_j}{S}-\frac{\alpha_k}{S})$-distributed.
$Beta\left(\frac{\alpha_k}{S},\sum_{j\neq k}\frac{\alpha_j}{S}\right)$-distributed.
\end{Rem}

\pff
We give here a short autonomous proof.
Denote $\alpha =\sum _{k=1}^d\alpha _k\geq 1$.
Conditional expectation at time $n+1$ writes
$$
\g E\left( P_{n+1}\left|\rond F_n\right.\right) =\frac{\alpha +(n+1)S}{\alpha +nS}P_n
$$
so that $\left(\frac{P_n}{\alpha +nS}\right) _{n\geq 0}$ is a $\left[ 0,1\right] ^d$-valued convergent martingale with mean $\left( \alpha _1/\alpha ,\dots ,\alpha _d/\alpha\right)$;
let $V$ be its limit.
If $f$ is any function defined on $\g R^d$,
$$
\g E\left( f\left( P_{n+1}\right)\left|\rond F_n\right.\right)
=\left( I+\frac{\Phi}{\alpha +nS}\right)\left( f\right)\left( P_n\right)
$$
where
$$
\Phi (f)(v)=\sum _{k=1}^{d}v_k\big[ f\left( v+Se_k\right) -f(v)\big]
$$
($e_k$ is the $k$-th vector in $\g R^{d}$ canonical basis and $v=\sum _{k=1}^dv_ke_k$).
In particular, as can be straightforwardly checked, if $p=\left( p_1,\dots ,p_d\right)\in\g N^d$ and $|p|=\sum _{k=1}^dp_k$, the function
$$
\Gamma _p(v)=
\prod _{k=1}^{d}\frac{\Gamma\left(\frac{v_k}{S}+p_k\right)}{\Gamma\left(\frac{v_k}{S}\right)},
$$
defined on $\g R^d$, is an eigenfunction of the operator $\Phi$, associated with the eigenvalue $|p|S$.
Consequently, after a direct induction, for any $p\in\g N^d$,
$$
\g E\left( \Gamma _p(P_n)\right)
=\frac{\Gamma\left(\frac \alpha S+n+|p|\right)}{\Gamma\left(\frac \alpha S+n\right)}
\cdot
\frac{\Gamma\left(\frac \alpha S\right)}{\Gamma\left(\frac \alpha S+|p|\right)}
\cdot
\Gamma _p(P_0)
%\prod _{k=1}^{S+1}\frac{\Gamma\left(\frac{1}{S}+p_k\right)}{\Gamma\left(\frac{1}{S}\right)} .
$$
so that, when $n$ tends to infinity, by Stirling's formula,
$$
\g E\left( \Gamma _p(P_n)\right)
=n^{|p]}
\cdot
\frac{\Gamma\left(\frac \alpha S\right)}{\Gamma\left(\frac \alpha S+|p|\right)}
\cdot
\Gamma _p\left( P_0\right)
\cdot
\left( 1+O\left( \frac 1n\right)\right) .
%\prod _{k=1}^{S+1}\frac{\Gamma\left(\frac{1}{S}+p_k\right)}{\Gamma\left(\frac{1}{S}\right)} .
$$
Besides, expanding real polynomials $X^p=X_1^{p_1}\dots X_{d}^{p_{d}}$ in the basis
$(\Gamma _p)_{p\in\g N^{d}}$, one gets formulae
$$
X^p=S^{|p|}\Gamma _p+\summ{k\in\g N^{d}}{|k|\leq |p|-1}a_{p,k}\Gamma _k(X)
$$
where the $a_{p,k}$ are rational numbers.
Consequently, when $n$ tends to infinity, one gets the asymptotics
$$
\g E\left( \frac{J_n}{\alpha +nS} \right)^p
=\frac{\Gamma\left(\frac \alpha S\right)}{\Gamma\left(\frac \alpha S+|p|\right)}
%\prod _{k=1}^{S+1}\frac{\Gamma\left(\frac{1}{S}+p_k\right)}{\Gamma\left(\frac{1}{S}\right)}
\Gamma _p(P_0)
\left( 1+O\left( \frac 1n\right)\right).
$$
which implies that, for any $p\in\g N^d$,
\begin{equation}
\label{momentsLimite}
\g E\left( V^p\right)
=\frac{\Gamma\left(\frac \alpha S\right)}{\Gamma\left(\frac \alpha S+|p|\right)}
\prod _{k=1}^{d}\frac{\Gamma\left(\frac{\alpha _k}{S}+p_k\right)}{\Gamma\left(\frac{\alpha _k}{S}\right)} .
\end{equation}
Note that this proves the convergence of the martingale in ${\rm L}^t$ for all $t\geq 1$.
Since a Dirichlet distribution is determined by its moments, this shows that the law of $V$ is a
Dirichlet distribution with parameters $\left(\frac{\alpha _1}{S},\dots ,\frac{\alpha _d}{S}\right)$.
\QED

{\bf Acknowledgements}

The authors wish to thank Quansheng Liu for stimulating discussions which gave birth to the recursive computation on the moments. Besides, one may have recognized Philippe Flajolet's style in the 
Maple figures which are originally due to him.

%%%%%%%%%%%%%%%%%%%%%%%%%%%%%%%%%%%%%%%%%%%%%
\bibliographystyle{plain}
\bibliography{bibCMP}
%%%%%%%%%%%%%%%%%%%%%%%%%%%%%%%%%%%%%%%%%%%%

\end{document}